\documentclass[11pt,letterpaper,colorlinks=true,linkcolor=blue,citecolor=blue,urlcolor=blue]{scrartcl}
\usepackage[T1]{fontenc}
\usepackage[utf8]{inputenc}
\usepackage{lmodern}
\usepackage[english]{babel}
\usepackage[english]{isodate}
\usepackage{amsthm}
\usepackage{amsmath}
\usepackage{amssymb}
\usepackage{amsfonts}
\usepackage{mathtools}
\usepackage{mathrsfs}
\usepackage{dsfont}
\usepackage{upgreek}
\usepackage{xpatch}
\usepackage{bm}
\usepackage{bbm}
\usepackage{enumitem}
\usepackage{float}
\usepackage[toc,page]{appendix}
\usepackage{footnotebackref}
\usepackage[totalwidth=480pt, totalheight=680pt]{geometry}
\usepackage{hyperref}
\usepackage{graphicx} 
\usepackage[font=footnotesize,skip=-18pt]{caption}

\theoremstyle{remark}
\newtheorem{remark}{Remark}[section]

\theoremstyle{definition}
\newtheorem{assumptions}[remark]{Assumptions}
\newtheorem{example}[remark]{Example}
\newtheorem{definition}[remark]{Definition}

\theoremstyle{plain}
\newtheorem{theorem}[remark]{Theorem}
\newtheorem{proposition}[remark]{Proposition}
\newtheorem{lemma}[remark]{Lemma}
\newtheorem{corollary}[remark]{Corollary}

\addtokomafont{disposition}{\normalfont\scshape}

\makeatletter
\xpatchcmd{\@sec@pppage}{
\bfseries}{
\normalfont\scshape\Large}{}{}
\makeatother

\setkomafont{section}{\large}
\setkomafont{subsection}{\normalsize}

\numberwithin{equation}{section}

\newcommand{\mybar}[1]{\makebox[0pt]{$\phantom{#1}\overline{\phantom{#1}}$}#1}
\newcommand{\myeqref}[1]{\hyperref[#1]{(\ref*{#1})}}
\newcommand{\myassum}[2]{\hyperref[#2]{Assumptions \ref*{#1} \ref*{#2}}}
\newcommand{\myassumption}[1]{\hyperref[#1]{Assumptions \ref*{#1}}}
\newcommand{\myappref}[1]{\hyperref[#1]{Appendix \ref*{#1}}}
\newcommand{\mysecref}[1]{\hyperref[#1]{Section \ref*{#1}}}
\newcommand{\mysubsecref}[1]{\hyperref[#1]{Subsection \ref*{#1}}}
\newcommand{\mythmref}[1]{\hyperref[#1]{Theorem \ref*{#1}}}
\newcommand{\mycorref}[1]{\hyperref[#1]{Corollary \ref*{#1}}}
\newcommand{\mylemref}[1]{\hyperref[#1]{Lemma \ref*{#1}}}
\newcommand{\myremref}[1]{\hyperref[#1]{Remark \ref*{#1}}}

\newcommand{\E}{\mathds{E}}
\newcommand{\W}{\mathds{W}}
\newcommand{\R}{\mathds{R}}
\newcommand{\G}{\mathds{G}}
\newcommand{\F}{\mathds{F}}
\renewcommand{\P}{\mathds{P}}

\renewcommand{\Pr}{\mathrm{P}}

\newcommand{\Ps}{\mathscr{P}}
\newcommand{\Fs}{\mathscr{F}}
\newcommand{\Ds}{\mathscr{D}}
\newcommand{\Bs}{\mathscr{B}}

\newcommand{\Uc}{\mathcal{U}}
\newcommand{\Vc}{\mathcal{V}}
\newcommand{\Wc}{\mathcal{W}}

\newcommand{\dd}{\mathrm{d}}
\newcommand{\bbeta}{\boldsymbol{\beta}}

\begin{document}


\title{
\LARGE A trajectorial approach to relative entropy dissipation of McKean--Vlasov diffusions: gradient flows and HWBI inequalities}


\author{
\large Bertram Tschiderer\thanks{
Faculty of Mathematics, University of Vienna, Oskar-Morgenstern-Platz 1, 1090 Vienna, Austria 
\newline (email: \href{mailto: bertram.tschiderer@univie.ac.at}{bertram.tschiderer@univie.ac.at}).
} 
\and \large Lane Chun Yeung\thanks{
Department of Industrial Engineering and Operations Research, Columbia University, 2990 Broadway, New York, NY 10027, USA 
\newline (email: \href{mailto: l.yeung@columbia.edu}{l.yeung@columbia.edu}).
} 
}


\date{\normalsize 26th May 2021}


\maketitle


\begin{abstract} \small \noindent \textsc{Abstract.} We formulate a trajectorial version of the relative entropy dissipation identity for McKean--Vlasov diffusions, extending the results of the papers \cite{FJ16,KST20a}, which apply to non-interacting diffusions. Our stochastic analysis approach is based on time-reversal of diffusions and Lions' differential calculus over Wasserstein space. It allows us to compute explicitly the rate of relative entropy dissipation along every trajectory of the underlying diffusion via the semimartingale decomposition of the corresponding relative entropy process. As a first application, we obtain a new interpretation of the gradient flow structure for the granular media equation, generalizing the formulation in \cite{KST20a} developed for the linear Fokker--Planck equation. Secondly, we show how the trajectorial approach leads to a new derivation of the HWBI inequality, which relates relative entropy (H), Wasserstein distance (W), barycenter (B) and Fisher information (I).

\bigskip

\small \noindent \href{https://mathscinet.ams.org/mathscinet/msc/msc2020.html}{\textit{MSC 2020 subject classifications:}} Primary 60H30, 60G44; secondary 60J60, 94A17

\bigskip

\small \noindent \textit{Keywords and phrases:} relative entropy dissipation, gradient flow, McKean--Vlasov diffusion, granular media equation, HWBI inequality
\end{abstract}


\section{Introduction}


We are interested in the relative entropy dissipation of McKean--Vlasov stochastic differential equations of the form 
\begin{equation} \label{eq: intro.1.1}
\mathrm{d}X_{t} = - \big( \nabla V(X_{t}) + \nabla (W \ast \Pr_{\kern-0.1em t})(X_{t}) \big) \, \mathrm{d}t + \sqrt{2} \, \mathrm{d}B_{t} \, , \qquad 0 \leqslant t \leqslant T,
\end{equation}
where $X_{0}$ has some given initial distribution $\Pr_{\kern-0.1em 0}$ on $\mathds{R}^{n}$. Here, the functions $V,W \colon \mathds{R}^{n} \rightarrow [0,\infty)$ play the roles of confinement and interaction potentials and are assumed to be suitably regular, $\Pr_{\kern-0.1em t} \coloneqq \operatorname{Law}(X_t)$ denotes the distribution of the random vector $X_{t}$, the symbol $\ast$ stands for the standard convolution operator, and $(B_{t})_{0 \leqslant t \leqslant T}$ is a standard $n$-dimensional Brownian motion. In particular, this SDE is non-local (or non-linear) in the sense that the drift term depends on the distribution of the state variable. Non-local equations of this form arise in the modeling of weakly interacting diffusion equations, after the seminal work of McKean \cite{McK66}. 

Since the work of Carrillo--McCann--Villani \cite{CMV03,CMV06}, relative entropy dissipation has been known to be an effective method for studying convergence rates to equilibrium and propagation of chaos of McKean--Vlasov equations. Some notable examples include the works \cite{CJMTU01,Mal03, BGG13, Tug13, CG14, CGPS20}. In a broader context, \cite{MMN18,HRSS20} recently applied entropy methods to the mean-field theory of neural networks.
 
We denote by $\mathscr{P}_{\mathrm{ac}}(\mathds{R}^{n})$ the set of absolutely continuous probability measures on $\mathds{R}^{n}$, which we will often identify with their corresponding probability density functions with respect to Lebesgue measure. The \emph{free energy functional} 
\begin{equation} \label{eq: free energy functional}
\mathscr{P}_{\mathrm{ac}}(\mathds{R}^{n}) \ni p \longmapsto \mathscr{F}(p) \coloneqq \Uc(p) + \Vc(p) + \Wc(p)
\end{equation} 
is defined as the sum of the energy functionals
\begin{equation} \label{eq: intro.1.4}
\mathcal{U}(p) \coloneqq \int_{\mathds{R}^{n}} p(x) \log p(x) \, \mathrm{d}x, \quad 
\mathcal{V}(p) \coloneqq \int_{\mathds{R}^{n}} V(x) \, p(x) \, \mathrm{d}x, \quad
\mathcal{W}(p) \coloneqq \tfrac{1}{2} \int_{\mathds{R}^{n}} (W\ast p)(x) \, p(x) \, \mathrm{d}x
\end{equation}
corresponding to internal ($\mathcal{U}$), potential ($\mathcal{V}$) and interaction ($\mathcal{W}$) energy, respectively. Defining the \emph{relative entropy dissipation functional} 
\begin{equation} \label{eq: entropy dissipation functional}
\mathscr{P}_{\mathrm{ac}}(\mathds{R}^{n}) \ni p \longmapsto \Ds(p) \coloneqq  \int_{\mathds{R}^{n}} \vert \nabla \log p(x) + \nabla V(x) + \nabla(W \ast p)(x) \vert^{2} \, p(x) \, \mathrm{d}x,
\end{equation}
the well-known relative entropy dissipation identity takes the form
\begin{equation} \label{eq: entropy dissipation identity}
\mathscr{F}(p_{t}) - \mathscr{F}(p_{t_{0}}) = - \int_{t_{0}}^{t}  \Ds(p_{u}) \, \mathrm{d}u.
\end{equation}
This identity is of a deterministic nature: it only depends on the curve of probability density functions $(p_{t})_{0 \leqslant t \leqslant T}$, but not on the trajectories of the underlying process $(X_{t})_{0 \leqslant t \leqslant T}$ itself. It is then natural to ask whether there is a process-level analogue of the relative entropy dissipation identity \hyperref[eq: entropy dissipation identity]{(\ref*{eq: entropy dissipation identity})}, depending directly on the trajectories of the McKean--Vlasov process $(X_{t})_{0 \leqslant t \leqslant T}$. The main contribution of this paper is to give an affirmative answer to this question, by formulating a trajectorial version of the relative entropy dissipation identity via a stochastic analysis approach.

Before going into details, let us briefly describe the main ideas. We draw inspiration from prior literature \cite{FJ16,KST20a} based on a simpler (linear) setting without interaction, i.e., $W \equiv 0$. In this case, the McKean--Vlasov SDE \hyperref[eq: intro.1.1]{(\ref*{eq: intro.1.1})} reduces to a Langevin--Smoluchowski diffusion equation of the form 
\begin{equation} \label{eq: l.s.e}
\mathrm{d}X_{t} = - \nabla V(X_{t}) \, \mathrm{d}t + \sqrt{2} \, \mathrm{d}B_{t} \, , \qquad 0 \leqslant t \leqslant T.
\end{equation}
In particular, the drift term does not depend on the distribution of $X_{t}$. Moreover, there is an explicit stationary distribution (also known as the Gibbs distribution) with density proportional to $q(x) \coloneqq \mathrm{e}^{-V(x)}$. Defining the likelihood ratio function (or Radon--Nikodym derivative) $\ell_t(x) \coloneqq p_t(x) \slash q(x)$, the free energy at time $t$ can be expressed as $\mathscr{F}(p_t) = \E[\log \ell_{t}(X_{t})]$, and the resulting stochastic process
\begin{equation}  \label{eq: r.e.p.}
\log \ell_{t}(X_{t}) = \log p_{t}(X_{t}) + V(X_{t}) \, , \qquad 0 \leqslant t \leqslant T
\end{equation}
is called \emph{free energy} or \emph{relative entropy process}. As shown in \cite{FJ16,KST20a}, the time-reversal $\log \ell_{T-s}(X_{T-s})_{0 \leqslant s \leqslant T}$ of this process is a submartingale, and It{\^o} calculus can be used to obtain its Doob--Meyer decomposition
\begin{equation} 
\log \ell_{T-s}(X_{T-s}) - \log \ell_{T}(X_{T}) = M_{T-s} + F_{T-s}.
\end{equation}
Here, $(M_{T-s})_{0 \leqslant s \leqslant T}$ is a martingale and $(F_{T-s})_{0 \leqslant s \leqslant T}$ is an increasing process of finite first variation, both with explicit expressions. This decomposition describes exactly the rate of relative entropy dissipation along every trajectory of the Langevin--Smoluchowski diffusion. Therefore, it can be viewed as a trajectorial analogue of the (deterministic) relative entropy dissipation identity \myeqref{eq: entropy dissipation identity}.

Let us now return to our McKean--Vlasov setting. In order to take into account the interaction potential $W$, it is natural to consider a generalized relative entropy process of the form
\begin{equation} \label{eq: e.r.e.p.}
\log p_{t}(X_{t}) + V(X_{t}) + \tfrac{1}{2} (W \ast \Pr_{\kern-0.1em t})(X_{t}) \, , \qquad 0 \leqslant t \leqslant T.
\end{equation}
The task is now to compute the semimartingale decomposition of this process. We will provide a detailed analysis of this extension, which is subtler than might appear at first sight. The main difficulty is that, even when it exists, the stationary distribution of the McKean--Vlasov diffusion does not have a closed-form expression in general. This prevents us from defining the likelihood ratio function in a straightforward manner as in the setting of Langevin--Smoluchowski diffusions, where one can rely on the invariant Gibbs distribution. An appropriate definition of the generalized likelihood ratio function turns out to be that \hyperref[eq: e.r.e.p.]{(\ref*{eq: e.r.e.p.})} should be viewed as a function of the form $\log \ell_{t}(X_{t}, \Pr_{\kern-0.1em t})$, depending explicitly on the distribution $\Pr_{\kern-0.1em t}$ of $X_{t}$ itself, in addition to the state $X_{t}$. This form of generalized likelihood ratio function allows us to take the L\emph{-derivative} with respect to the probability distribution $\Pr_{\kern-0.1em t}$. The notion of L-differentiation for functions of probability measures was introduced by Lions \cite{Lio08}. We refer to the monograph \cite[Chapter 5]{CD18a} for a detailed discussion of differential calculus and stochastic analysis over spaces of probability measures. In particular, we will use a generalized form of It\^{o}'s formula for functions of curves of measures, to derive the dynamics of the time-reversal of the relative entropy process \hyperref[eq: e.r.e.p.]{(\ref*{eq: e.r.e.p.})}, in terms of the semimartingale decomposition
\begin{equation} \label{eq: semi.decom}
\log \ell_{T-s}(X_{T-s}, \mathrm{P}_{\kern-0.1em T - s}) - \log \ell_{T}(X_{T}, \mathrm{P}_{\kern-0.1em T}) = M_{T-s} + F_{T-s} \, , \qquad 0 \leqslant s \leqslant T,
\end{equation}
where $(M_{T-s})_{0 \leqslant s \leqslant T}$ is a martingale and $(F_{T-s})_{0 \leqslant s \leqslant T}$ is a process of finite first variation, both of which will be explicitly computed. Similar to the case of Langevin--Smoluchowski dynamics, this decomposition can be viewed as the trajectorial rate of relative entropy dissipation. The classical (deterministic) identity \myeqref{eq: entropy dissipation identity} can then be recovered by taking expectations.


\subsection{Gradient flow structure of the granular media equation}


As a first application of our trajectorial approach we obtain a new interpretation of the gradient flow structure of the \emph{granular media equation}
\begin{equation} \label{eq: intro.1.2}
\partial_{t} p_{t}(x) = \operatorname{div} \Big( \nabla p_{t}(x) + p_{t}(x)  \nabla V(x) + p_{t}(x)  \nabla(W \ast p_{t})(x) \Big) \, , \qquad (t,x) \in (0,T) \times \mathds{R}^{n},
\end{equation}
which describes the evolution of the curve of probability density functions $(p_t)_{0 \leqslant t \leqslant T}$ corresponding to the McKean--Vlasov diffusion $(X_{t})_{0 \leqslant t \leqslant T}$ of \hyperref[eq: intro.1.1]{(\ref*{eq: intro.1.1})}. When $n=1$, this PDE appears in the modeling of the time evolution of granular media \cite{BCCP98,Vil06, CGM08}; in that context, the granular medium is modeled as system of particles performing inelastic collisions, and $p_{t}(x)$ is regarded as the velocity of a representative particle in the system at time $t$ and position $x$, while $V$ and $W$ represent the friction and the inelastic collision forces, respectively. Note that in the interaction-free case $W \equiv 0$, the equation \hyperref[eq: intro.1.2]{(\ref*{eq: intro.1.2})} reduces to a linear \emph{Fokker--Planck equation}. As is well known from \cite{CMV03,CMV06}, this curve of probability densities can be characterized as a \emph{gradient flow} in $\mathscr{P}_{\mathrm{ac},2}(\mathds{R}^{n})$, the space of absolutely continuous probability measures with finite second moments. Roughly speaking, this is an optimality property stating that the curve $(p_{t})_{0 \leqslant t \leqslant T}$ evolves in the direction of steepest possible descent for the free energy functional \hyperref[eq: free energy functional]{(\ref*{eq: free energy functional})} with respect to the quadratic Wasserstein distance 
\begin{equation} \label{eq: wasserstein dist}
W_{2}(\mu,\nu) = \Big( \, \inf_{\scriptscriptstyle Y \sim \mu, Z \sim \nu} \mathds{E} \vert Y - Z \vert^{2} \, \Big)^{1/2} \, , \qquad \mu, \nu \in \mathscr{P}_{2}(\mathds{R}^{n}).
\end{equation}

The Wasserstein gradient flow structure of the linear Fokker--Planck equation was first discovered by Jordan, Kinderlehrer and Otto in the seminal work \cite{JKO98}. In the paper \cite{OV00}, Otto and Villani developed a formal Riemannian structure on the space of probability measures with finite second moments, leading to heuristic proofs of gradient flow properties as in \cite{Ott01}, where the porous medium equation was studied. This pioneering approach is often referred to as ``Otto calculus''. Later, a rigorous framework based on minimizing movement schemes and curves of maximal slope was introduced in \cite{AGS08}. Recently, a trajectorial approach to the gradient flow properties of Langevin--Smoluchowski diffusions \cite{KST20a} and Markov chains \cite{KMS20} was established. We will follow this approach and adapt it to our McKean--Vlasov setting. For gradient flows of McKean--Vlasov equations on discrete spaces we refer to \cite{EFLS16}.

Returning to the setting of this paper, our main result leads to a new formulation of the gradient flow property of the granular media equation. To show this steepest descent property, the main idea is to consider a \emph{perturbed} McKean--Vlasov diffusion of the form
\begin{equation} \label{eq: intro.p.sde}
\mathrm{d}X_{t} = - \big( \nabla V(X_{t}) + \nabla \bbeta(X_{t}) + \nabla (W \ast \mathrm{P}_{\kern-0.1em t}^{\boldsymbol{\beta}})(X_{t}) \big) \, \mathrm{d}t + \sqrt{2} \, \mathrm{d}B_{t}^{\boldsymbol{\beta}}
\, , \qquad t_{0} \leqslant t \leqslant T
\end{equation}
which is constructed by adding a perturbation $\bbeta \colon \R^{n} \rightarrow \R$ to the confinement potential\footnote{As we will see, the steepest descent property is already visible by perturbing the confinement potential from $V$ to $V+\bbeta$, thus we avoid complicating the setup further by adding another perturbation to the interaction potential $W$.} of the original McKean--Vlasov SDE \myeqref{eq: intro.1.1}. In other words, from time $t_0$ onward, the perturbed diffusion drifts in a direction different from that of the original diffusion, hence the perturbed curve of time-marginal distributions $(\Pr^{\bbeta}_{\kern-0.1em t})_{t_{0} \leqslant t \leqslant T}$ also evolves differently from the unperturbed curve $(\Pr_{\kern-0.1em t})_{t_{0} \leqslant t \leqslant T}$. In parallel with the unperturbed case, we may compute the dynamics of the perturbed relative entropy process associated with \myeqref{eq: intro.p.sde}. As a consequence, we derive the rate of relative entropy dissipation for the perturbed McKean--Vlasov diffusion. On the other hand, the rate of change of the Wasserstein distance along the perturbed curve $(\Pr^{\bbeta}_{\kern-0.1em t})_{t_{0} \leqslant t \leqslant T}$ can be computed based on the general theory of metric derivative of absolutely continuous curves, see \cite{AGS08}. Finally, comparing these two rates in both the perturbed and unperturbed settings, allows us to establish the gradient flow property.


\subsection{The HWBI inequality}


The second application of our trajectorial approach deals with the HWBI inequality \cite[Theorem 4.2]{AGK04}, which is an extension of the HWI inequality \cite{OV00}. It relates not only relative entropy (H), Wasserstein distance (W), and relative Fisher information (I), but also barycenter (B). These quantities are defined as follows: for two probability measures $\nu, \mu \in \mathscr{P}(\mathds{R}^{n})$, the \emph{relative entropy} of $\nu$ with respect to $\mu$ is defined by
\begin{equation} \label{eq: def.rel.ent}
H(\nu \, \vert \, \mu) \coloneqq 
\begin{cases}
\int_{\mathds{R}^{n}} \frac{\mathrm{d}\nu}{\mathrm{d}\mu} \log (\frac{\mathrm{d}\nu}{\mathrm{d}\mu}) \, \mathrm{d}\mu, \, & \text{if } \nu \ll \mu \\
+ \infty, & \text{otherwise},
\end{cases}
\end{equation}
the \emph{relative Fisher information} of $\nu$ with respect to $\mu$ is given by
\begin{equation} \label{eq: def.rel.fish}
I(\nu \, \vert \, \mu) \coloneqq 
\begin{cases}
\int_{\mathds{R}^{n}} \vert \nabla \log (\frac{\mathrm{d}\nu}{\mathrm{d}\mu})\vert^{2}  \, \mathrm{d}\mu, \,  & \text{if }\nu \ll \mu \\
+ \infty, & \text{otherwise},
\end{cases}
\end{equation}
and the \emph{barycenter} of a probability measure $\nu \in \mathscr{P}_{2}(\mathds{R}^{n})$ is defined as $b(\nu) \coloneqq \int_{\mathds{R}^{n}} x  \, \mathrm{d}\nu(x) \in \mathds{R}^{n}$, where the integral is understood as a Bochner integral. Informally, the HWBI inequality then states that any two probability measures $\nu_{0}, \nu_{1} \in \mathscr{P}_{2}(\mathds{R}^{n})$ satisfy
\begin{equation} \label{eq: HWBI.o}
H(\nu_{0} \, \vert \, \mu_{0}) - H(\nu_{1} \, \vert \, \mu_{1}) \leqslant
\sqrt{I(\nu_{0} \, \vert \, \mu_{0}^{\uparrow})} \, W_{2}(\nu_{0},\nu_{1})
- \tfrac{\kappa_{V} + \kappa_{W}}{2} \, W_{2}^{2}(\nu_{0},\nu_{1}) + \tfrac{\kappa_{W}}{2} \, \vert b(\nu_{0}) - b(\nu_{1})\vert^{2},
\end{equation}
where $\mu_{0}, \mu_{1}, \mu_{0}^{\uparrow}$ are some appropriate $\sigma$-finite reference measures depending on the potentials $V,W$ (see \hyperref[subs: hwbi]{Subsection \ref*{subs: hwbi}} for the details), and $\kappa_{V}, \kappa_{W} \in \mathds{R}$ are the moduli of uniform convexity for $V,W$. This inequality describes the evolution of the relative entropy along the \emph{displacement interpolation}  $(\nu_t)_{0 \leqslant t \leqslant 1}$ between $\nu_{0}$ and $\nu_{1}$. Compared with the HWI inequality, there are two additional terms on the right-hand side of \myeqref{eq: HWBI.o} contributed by the interaction energy  functional $\mathcal{W}$ of \myeqref{eq: intro.1.4}. Intuitively, the $\kappa_W$-uniform convexity of $W$ leads to the first additional term $-\tfrac{\kappa_W}{2} W_2^2(\nu_0, \nu_1)$, which alone would correspond to the $\kappa_W$-uniform displacement convexity of $\mathcal{W}$ along $(\nu_t)_{0 \leqslant t \leqslant 1}$. But since $\mathcal{W}(p)$ is invariant under any translation of $p$, the functional $\mathcal{W}$ might fail to be uniformly displacement convex when the barycenter shifts. This suggests that the barycentric shift along $(\nu_t)_{0 \leqslant t \leqslant 1}$ should be factored out of the consideration of the displacement convexity of $\mathcal{W}$, which is intuitively why the second additional term $\frac{\kappa_W}{2} \vert b(\nu_{0}) - b(\nu_{1})\vert^{2}$ in \myeqref{eq: HWBI.o} appears. 

Coming back to our second application, we illustrate how our approach yields a trajectorial proof of the inequality \hyperref[eq: HWBI.o]{(\ref*{eq: HWBI.o})}, in the slightly strengthened form of \cite[Theorem 4.1]{CEGH04} and \cite[Theorem D.50]{FK06}. Much of this consists of arguments similar in spirit to our main result \hyperref[eq: semi.decom]{(\ref*{eq: semi.decom})}, but with one key difference: instead of the time-marginals of the McKean--Vlasov diffusion, we apply the trajectorial approach to the displacement interpolation $(\nu_t)_{0 \leqslant t \leqslant 1}$. In this regard, our derivation can be seen as a generalization of the trajectorial proof of the HWI inequality in \cite[Section 4.2]{KST20a}; see also \cite[Section 9.4]{KMS20}, where the same idea was used to derive a discrete version of the HWI inequality in a Riemannian-geometric framework.

In the literature, similar trajectorial approaches have also been applied in the context of martingale inequalities \cite{ABPST13,BS15}, functional inequalities \cite{Cat04,Leh13,BCGL20,GLRT20}, and their stability estimates \cite{ELS20,EM20}. In particular, we refer to \cite[Corollary 1.4]{BCGL20} for a related HWI inequality derived from the entropic interpolation of the mean-field Schr\"odinger problem. 


\subsection{Organization of the paper}


We set up the probabilistic framework and discuss some regularity assumptions in \hyperref[sec: framework]{Section \ref*{sec: framework}}. In \mysecref{sec: main results} we state our main trajectorial results, \mythmref{thm: trac} and \mythmref{thm: p.trac}, and develop two explicit examples for illustration. As immediate consequences, we derive the classical relative entropy dissipation identities in \mycorref{cor: ent.diss.id} and \mycorref{cor: per.rel.ent.id.b}. Building on these results, we formulate the gradient flow property of the granular media equation in \mythmref{prop: gradient flow}. The HWBI inequality is then stated in \mythmref{thm: HWBI}. The proofs of the trajectorial results and of the HWBI inequality are developed in \hyperref[sec: Proofs]{Section \ref*{sec: Proofs}}. Some proofs of auxiliary results postponed in previous sections are contained in \hyperref[app]{Section \ref*{app}}.


\section{The probabilistic framework} \label{sec: framework}


\subsection{The setting} \label{sec: notations}


We fix a terminal time $T \in (0,\infty)$ and let $\Omega \coloneqq C([0,T];\mathds{R}^{n})$ be the path space of $\mathds{R}^{n}$-valued continuous functions defined on $[0,T]$. We denote by $(X_{t})_{0 \leqslant t \leqslant T}$ the canonical process defined by $X_{t}(\omega) \coloneqq \omega(t)$ for $\omega \in \Omega$, and fix a probability distribution $\mathrm{P}_{\kern-0.1em 0} \in \mathscr{P}_{\mathrm{ac},2}(\mathds{R}^{n})$.

As will be shown in \hyperref[lem: well-posedness]{Lemma \ref*{lem: well-posedness}}, under the \textnormal{\hyperref[assumptions]{Assumptions \ref*{assumptions}}} below, the SDE \hyperref[eq: intro.1.1]{(\ref*{eq: intro.1.1})} with initial distribution $\mathrm{P}_{\kern-0.1em 0}$ has a unique strong solution, when it is posed on an arbitrary filtered probability space. This implies that there exists a probability measure $\P$ on $\Omega$ and a $\P$-Brownian motion $(B_{t})_{0 \leqslant t \leqslant T}$ such that the SDE \hyperref[eq: intro.1.1]{(\ref*{eq: intro.1.1})} holds. We write $\F = (\mathcal{F}_{t})_{0 \leqslant t \leqslant T}$ for the right-continuous augmentation of the canonical filtration. 

For each time $t \in [0,T]$, we denote by $\mathrm{P}_{\kern-0.1em t} \coloneqq \P \circ X_{t}^{-1}$ the distribution of $X_{t}$ under $\P$, and by $p_{t}$ the corresponding probability density function on $\mathds{R}^{n}$. The density functions $(p_{t})_{0 \leqslant t \leqslant T}$ then solve the granular media equation \hyperref[eq: intro.1.2]{(\ref*{eq: intro.1.2})}.  


\subsection{Regularity assumptions}  \label{subsec: regularity assumptions}


\begin{assumptions} \label{assumptions} 
The following regularity assumptions will be used frequently.
\begin{enumerate}[label=(\roman*)] 
\item \label{r.a.1} The functions $V,W \colon \mathds{R}^{n} \rightarrow [0,\infty)$ are smooth and have Lipschitz continuous gradients with Lipschitz constants $\Vert \nabla V \Vert_{\textnormal{Lip}}, \Vert \nabla W \Vert_{\textnormal{Lip}}$. All derivatives of $V$ and $W$ grow at most exponentially as $\vert x \vert$ tends to infinity, and the first derivatives are of linear growth. The latter condition means that there exists a constant $C > 0$ such that 
\begin{equation} \label{lin.gro.con.}
\forall \, x \in \mathds{R}^{n} \colon \qquad \vert \nabla V(x) \vert \leqslant C (1 + \vert x \vert) \, , \qquad \vert \nabla W(x) \vert \leqslant C (1 + \vert x \vert).
\end{equation}
Furthermore, the function $W$ is even (in other words, symmetric), i.e., $W(x) = W(-x)$ for all $x \in \mathds{R}^{n}$.
\item \label{r.a.4} The probability distribution $\Pr_{\kern-0.1em 0}$ is an element of the space $\mathscr{P}_{\kern-0.1em \mathrm{ac},2}(\mathds{R}^{n})$ and the corresponding probability density function $\mathds{R}^{n} \ni x \mapsto p_{0}(x)$ is strictly positive. Moreover, the initial free energy $\mathscr{F}(p_{0})$ is finite.
\end{enumerate}
\end{assumptions}

These assumptions ensure that the equation $\hyperref[eq: intro.1.1]{(\ref*{eq: intro.1.1})}$ belongs to a broad class of strongly solvable McKean--Vlasov SDEs. We relegate the proof of the following result to \hyperref[sec: pf.well-pose]{Subsection \ref*{sec: pf.well-pose}}. 

\begin{lemma} \label{lem: well-posedness} Suppose \textnormal{\hyperref[assumptions]{Assumptions \ref*{assumptions}}} hold. Then on an arbitrary filtered probability space, the McKean--Vlasov SDE \textnormal{\myeqref{eq: intro.1.1}} has a pathwise unique, strong solution $(X_{t})_{0 \leqslant t \leqslant T}$ satisfying
\begin{equation} \label{fsmc}
\E  \bigg[ \sup_{0 \leqslant t \leqslant T} \vert X_{t} \vert^{2}\bigg] < \infty.
\end{equation}
Moreover, its marginal distributions $(\mathrm{P}_{\kern-0.1em t})_{0 \leqslant t \leqslant T}$ belong to $\mathscr{P}_{\kern-0.1em \mathrm{ac},2}(\mathds{R}^{n})$, and the corresponding curve of probability density functions $(p_{t})_{0 \leqslant t \leqslant T}$ is a classical solution of the granular media equation \textnormal{\hyperref[eq: intro.1.2]{(\ref*{eq: intro.1.2})}}.
\end{lemma}


\subsection{Probabilistic representations of gradient flow functionals}


To set up our framework, the first step is to express the free energy as well as the relative entropy dissipation functional in probabilistic terms. To this end, we introduce the generalized potential $\Psi \colon \mathds{R}^{n} \times \mathscr{P}_{2}(\mathds{R}^{n}) 
\to [0, \infty)$ and its close relative $\Psi^{\uparrow}$ given by
\begin{equation} \label{eq: pot}
\Psi(x,\mu) \coloneqq V(x) + \tfrac{1}{2}(W\ast\mu)(x) \, , \qquad 
\Psi^{\uparrow}(x,\mu) \coloneqq V(x) + (W\ast\mu)(x) \, 
\end{equation}
for $(x,\mu) \in \mathds{R}^{n} \times \mathscr{P}_{2}(\mathds{R}^{n})$. Furthermore, we define the density functions
\begin{equation} \label{eq: d.f.s.}
q(x,\mu) \coloneqq \mathrm{e}^{-\Psi(x,\mu)} \ , \qquad \qquad
q^{\uparrow}(x,\mu) \coloneqq \mathrm{e}^{-\Psi^{\uparrow}(x,\mu)} \ , \qquad \qquad
q^{\downarrow}(x) \coloneqq \mathrm{e}^{-V(x)}
\end{equation}
and the corresponding generalized likelihood ratio functions
\begin{equation} \label{eq: intro.1.7}
\ell_{t}(x,\mu) \coloneqq \frac{p_{t}(x)}{q(x,\mu)} \ , \qquad \qquad
\ell_{t}^{\uparrow}(x,\mu) \coloneqq \frac{p_{t}(x)}{q^{\uparrow}(x,\mu)} \ , \qquad \qquad
\ell_{t}^{\downarrow}(x) \coloneqq \frac{p_{t}(x)}{q^{\downarrow}(x)}
\end{equation}
for $t \in [0,T]$. Note that if $W \equiv 0$, these three likelihood ratio functions coincide. 

For each time $t \in [0, T]$, we introduce $\sigma$-finite measures on the Borel sets of $\mathds{R}^{n}$, given by
\begin{equation} \label{eq: intro.1.6}
\mathrm{Q}_{t}(A) \coloneqq \int_{A} q(x,\mathrm{P}_{\kern-0.1em t}) \, \mathrm{d}x 
\, , \qquad  
\mathrm{Q}_{t}^{\uparrow}(A) \coloneqq \int_{A} q^{\uparrow}(x,\mathrm{P}_{\kern-0.1em t}) \, \mathrm{d}x 
\, , \qquad 
A \in \mathscr{B}(\mathds{R}^{n}).
\end{equation}
Intuitively, these measures are (unnormalized) time-dependent Gibbs distributions. If $W \equiv 0$, they coincide with the true Gibbs distribution of the Langevin--Smoluchowski equation \myeqref{eq: l.s.e}, which is also its stationary distribution (when normalized to a probability measure).

With these definitions, we can now write the gradient flow functionals $\mathscr{F}$ and $\mathscr{D}$, introduced in \myeqref{eq: free energy functional} and \myeqref{eq: entropy dissipation functional}, in probabilistic terms: the relative entropy (defined in \myeqref{eq: def.rel.ent}) of $\mathrm{P}_{\kern-0.1em t}$ with respect to $\mathrm{Q}_{t}$ and the relative Fisher information (defined in  \myeqref{eq: def.rel.fish}) of $\mathrm{P}_{\kern-0.1em t}$ with respect to $\mathrm{Q}_{t}^{\uparrow}$ can be expressed respectively as
\begin{equation} \label{eq: intro.1.8}
H(\mathrm{P}_{\kern-0.1em t} \, \vert \, \mathrm{Q}_{t}) = \mathds{E}_{\mathds{P}}\big[\log\ell_{t}(X_{t},\mathrm{P}_{\kern-0.1em t})\big]
\, , \qquad 
I(\mathrm{P}_{\kern-0.1em t} \, \vert \, \mathrm{Q}_{t}^{\uparrow}) 
= \mathds{E}_{\mathds{P}}\big[\vert\nabla\log\ell_{t}^{\uparrow}(X_{t},\mathrm{P}_{\kern-0.1em t})\vert^{2}\big] \, ;
\end{equation}
and we have the relations $H(\mathrm{P}_{\kern-0.1em t} \, \vert \, \mathrm{Q}_{t}) = \mathscr{F}(\mathrm{P}_{\kern-0.1em t})$ as well as $I(\mathrm{P}_{\kern-0.1em t} \, \vert \, \mathrm{Q}_{t}^{\uparrow}) = \Ds(\mathrm{P}_{\kern-0.1em t})$. In particular, the relative entropy $H(\mathrm{P}_{\kern-0.1em t} \, \vert \, \mathrm{Q}_{t})$ can be written as the $\mathds{P}$-expectation of the \emph{relative entropy process} 
\begin{equation} \label{eq: intro.2.11}
\log \ell_{t} (X_{t},\mathrm{P}_{\kern-0.1em t}) = \log p_{t}(X_{t}) + V(X_{t}) + \tfrac{1}{2}(W \ast \mathrm{P}_{\kern-0.1em t})(X_{t}) \, , \qquad 0 \leqslant t \leqslant T.
\end{equation}
The dynamics of this stochastic process, together with its perturbed counterpart to be introduced in \hyperref[subsec: gradient flow]{Subsection \ref*{subsec: gradient flow}} below, will be our main objects of interest.

\begin{remark} \label{rm: rel ent} If the reference measure $\mathrm{Q}_{t}$ in \hyperref[eq: intro.1.6]{(\ref*{eq: intro.1.6})} is a \emph{probability} measure, then the expression \hyperref[eq: intro.1.8]{(\ref*{eq: intro.1.8})} matches the classical definition of relative entropy given in \hyperref[eq: def.rel.ent]{(\ref*{eq: def.rel.ent})}. In the general case when $\mathrm{Q}_{t}$ is a $\sigma$-finite measure, the definition \hyperref[eq: def.rel.ent]{(\ref*{eq: def.rel.ent})} is also valid under the condition that $\mathrm{P}_{\kern-0.1em t}$ has finite second moment, with the only difference that the range of the function $t \mapsto H(\mathrm{P}_{\kern-0.1em t} \, \vert \, \mathrm{Q}_{t})$ is extended from $[0,\infty]$ to $(-\infty,\infty]$; we refer to \cite[Appendix C]{KST20b} or \cite[Section 3]{Leo14b} for the details. 
\end{remark}


\section{Main results} \label{sec: main results}


\subsection{Trajectorial dissipation of relative entropy for McKean--Vlasov diffusions} \label{subsec: traj diss}


Our first main result is the semimartingale decomposition of the relative entropy process \hyperref[eq: intro.2.11]{(\ref*{eq: intro.2.11})}. It describes the dissipation of relative entropy along every trajectory of a particle undergoing the McKean--Vlasov dynamics \hyperref[eq: intro.1.1]{(\ref*{eq: intro.1.1})}. In the same spirit as the trajectorial approaches of \cite{FJ16} and \cite{KST20a}, we shall study the dynamics of the relative entropy process in the backward direction of time. Concretely, we consider for arbitrary, fixed $T \in (0,\infty)$ the time-reversed canonical process
\begin{equation} \label{eq: t.r.c.p.}
\mybar{X}_{s} \coloneqq X_{T-s} \, , \qquad 0 \leqslant s \leqslant T
\end{equation}
on the filtered probability space $(\Omega,\mathds{G},\mathds{P})$, where $\mathds{G} = (\mathcal{G}_{s})_{0 \leqslant s \leqslant T}$ is the $\mathds{P}$-augmented filtration generated by $(\mybar{X}_{s})_{0 \leqslant s \leqslant T}$. 

In order to formulate \mythmref{thm: trac} below, we introduce the time-reversed \emph{Fisher information process} 
\begingroup
\addtolength{\jot}{0.7em} 
\begin{equation} \label{eq: f.i.p.}
\begin{aligned}
\mybar{I}_{s} 
&\coloneqq \Big( \big\vert \nabla\log\mybar{\ell}_{s}^{\downarrow} \big\vert^{2} + \tfrac{1}{2}  \big\vert \nabla( W\ast\mybar{\mathrm{P}}_{\kern-0.1em s}) \big\vert^{2} + \Big\langle  \tfrac{1}{2} \nabla(W\ast\mybar{\mathrm{P}}_{\kern-0.1em s}) \,  , \, 2\nabla\log\mybar{\ell}_{s}^{\downarrow} + \nabla V  \Big\rangle \Big)(\mybar{X}_{s})  \\ 
& \qquad -  \mathds{E}_{\tilde{\mathds{P}}}\Big[\Big\langle  \tfrac{1}{2} \nabla W(\mybar{X}_{s}-\mybar{Y}_{s})  \, , \Big(2\nabla\log\mybar{\ell}_{s}^{\downarrow} - \nabla V +\nabla ( W\ast\mybar{\mathrm{P}}_{\kern-0.1em s})\Big)(\mybar{Y}_{s})\Big\rangle \Big] 
\end{aligned}
\end{equation}
\endgroup
for $0 \leqslant s \leqslant T$. Here, the process $(\mybar{Y}_{s})_{0 \leqslant s \leqslant T}$ is defined on another probability space $(\tilde{\Omega},\tilde{\mathds{G}},\tilde{\mathds{P}})$ such that the tuple $(\tilde{\Omega},\tilde{\mathds{G}},\tilde{\mathds{P}},(\mybar{Y}_{s})_{0 \leqslant s \leqslant T})$ is an exact copy of $(\Omega,\mathds{G},\mathds{P},(\mybar{X}_{s})_{0 \leqslant s \leqslant T})$. A bar over a letter means that time is reversed as in \hyperref[eq: t.r.c.p.]{(\ref*{eq: t.r.c.p.})}.

We also define the time-reversed \emph{cumulative Fisher information process} as the time integral
\begin{equation} \label{eq: c.f.i.}
\mybar{F}_{s} \coloneqq  \int_{0}^{s} \mybar{I}_{u} \, \mathrm{d}u  \, , \qquad 0 \leqslant s \leqslant T.
\end{equation}
This process will act as the compensator in the semimartingale decomposition of the relative entropy process \hyperref[eq: intro.2.11]{(\ref*{eq: intro.2.11})}. Its relation with the relative Fisher information \hyperref[eq: intro.1.8]{(\ref*{eq: intro.1.8})} will be given in \hyperref[eq: fish.inf.exp]{(\ref*{eq: fish.inf.exp})} below.

\begin{theorem} \label{thm: trac} Suppose \textnormal{\myassumption{assumptions}} hold. On $(\Omega,\mathds{G},\mathds{P})$, the time-reversed relative entropy process 
\begin{equation} \label{eq: t.r.r.e.pr}
\log \mybar{\ell}_{s}(\mybar{X}_{s},\mybar{\mathrm{P}}_{\kern-0.1em s}) = \log\ell_{T-s}(X_{T-s},\mathrm{P}_{\kern-0.1em T-s}) \, , \qquad  0 \leqslant s \leqslant T
\end{equation}
admits the semimartingale decomposition 
\begin{equation} \label{eq: d.m.d}
\log \mybar{\ell}_{s}(\mybar{X}_{s},\mybar{\mathrm{P}}_{\kern-0.1em s}) - \log \mybar{\ell}_{0}(\mybar{X}_{0},\mybar{\mathrm{P}}_{\kern-0.1em 0}) = \mybar{M}_{s} + \mybar{F}_{s}.
\end{equation}
Here $(\mybar{M}_{s})_{0 \leqslant s \leqslant T}$ is the $L^{2}(\mathds{P})$-bounded martingale
\begin{equation} \label{eq: martingale}
\mybar{M}_{s} \coloneqq \int_{0}^{s} \Big\langle \nabla \log \mybar{\ell}_{u}(\mybar{X}_{u},\mybar{\mathrm{P}}_{\kern-0.1em u}) \, , \, \sqrt{2} \, \mathrm{d}\mybar{B}_{u} \Big\rangle,
\end{equation}
with $(\mybar{B}_{s})_{0 \leqslant s \leqslant T}$ a $\P$-Brownian motion of the backward filtration $\G$, and the compensator $(\mybar{F}_{s})_{0 \leqslant s \leqslant T}$ satisfies 
\begin{equation} \label{eq: fish.inf.exp}
\mathds{E}_{\mathds{P}} \big[\mybar{F}_{s}\big] 
= \int_{0}^{s} I\big(\mybar{\mathrm{P}}_{\kern-0.1em u} \, \big\vert \, \mybar{\mathrm{Q}}_{u}^{\uparrow}\big) \, \mathrm{d}u 
= \mathds{E}_{\mathds{P}}\bigg[ \int_{0}^{s} \big\vert\nabla\log\mybar{\ell}_{u}^{\uparrow}(\mybar{X}_{u},\mybar{\mathrm{P}}_{\kern-0.1em u})\big\vert^{2} \, \mathrm{d}u \bigg] < \infty.
\end{equation}
\end{theorem}


\subsubsection{Examples}


We give two concrete examples to illustrate \hyperref[thm: trac]{Theorem \ref*{thm: trac}}.

\begin{example} \label{exa: o.u.} We set $n = 1$ and specialize \hyperref[thm: trac]{Theorem \ref*{thm: trac}} to the case of quadratic confinement potential $V(x) = \frac{x^{2}}{2}$ and no interaction potential $W \equiv 0$. The initial position $X_{0}$ is chosen to be independent of $(B_{t})_{0 \leqslant t \leqslant T}$ and to be normally distributed with mean $0$ and variance $\sigma_{0}^{2} > 0$. In this case, the SDE of \hyperref[eq: intro.1.1]{(\ref*{eq: intro.1.1})} becomes
\begin{equation} \label{eq: o.u.}
\mathrm{d}X_{t} = - X_{t} \, \mathrm{d}t + \sqrt{2} \, \mathrm{d}B_{t} \, , \qquad 0 \leqslant t \leqslant T
\end{equation}
and its solution is given by the Ornstein--Uhlenbeck process 
\begin{equation} \label{eq: sol.ou.pr}
X_{t} = \mathrm{e}^{-t} X_{0} + \sqrt{2} \int_{0}^{t} \mathrm{e}^{u-t} \, \mathrm{d}B_{u} \, , \qquad 0 \leqslant t \leqslant T
\end{equation}
with probability density function 
\begin{equation} \label{eq: o.u.p.d.f.}
p_{t}(x) = \tfrac{1}{\sqrt{2\pi\sigma_{t}^{2}}} \exp\big(-\tfrac{x^{2}}{2\sigma_{t}^{2}}\big)
\, , \qquad 
\sigma_{t}^{2} \coloneqq  1 + \mathrm{e}^{-2t}(\sigma_{0}^{2} - 1).
\end{equation}
Recalling \hyperref[eq: pot]{(\ref*{eq: pot})} -- \hyperref[eq: intro.1.7]{(\ref*{eq: intro.1.7})} and using \hyperref[eq: o.u.p.d.f.]{(\ref*{eq: o.u.p.d.f.})}, we see that in this setting the cumulative Fisher information process of \hyperref[eq: c.f.i.]{(\ref*{eq: c.f.i.})} is explicitly given by 
\begin{equation} \label{eq: o.u.f.}
\mybar{F}_{s}^{\mathrm{OU}} = \int_{0}^{s}  \Big( \nabla\log\mybar{\ell}_{u}^{\downarrow}(\mybar{X}_{u}) \Big)^{2} \, \mathrm{d}u
= \int_{0}^{s}  \Big( \nabla\log\mybar{p}_{u}(\mybar{X}_{u}) + \mybar{X}_{u} \Big)^{2} \, \mathrm{d}u
= \int_{0}^{s}  \big(1-\tfrac{1}{\overline{\sigma}_{u}^{2}}\big)^{2} \, \mybar{X}_{u}^{2}  \, \mathrm{d}u
\end{equation}
for $0 \leqslant s \leqslant T$. In particular, the non-negativity of the integrand in \hyperref[eq: o.u.f.]{(\ref*{eq: o.u.f.})} implies that the relative entropy decreases along almost every trajectory. 

Now we set $T = 1$ and $\sigma_{0}^{2} = 0.1$. The blue lines in \hyperref[fig: OU.sim]{Figure \ref*{fig: OU.sim}} below are ten simulated trajectories $s \mapsto \mybar{F}_{s}^{\mathrm{OU}}(\omega_{i})$, for $i = 1, \ldots, 10$. The thick black line plots the expected path $s \mapsto \mathds{E}_{\mathds{P}}[\mybar{F}_{s}^{\mathrm{OU}}]$ of all possible trajectories. 

\vspace*{-1cm}

\begin{figure}[H] 
\centering
\captionsetup{justification=centering}
\includegraphics[width=0.6\textwidth]{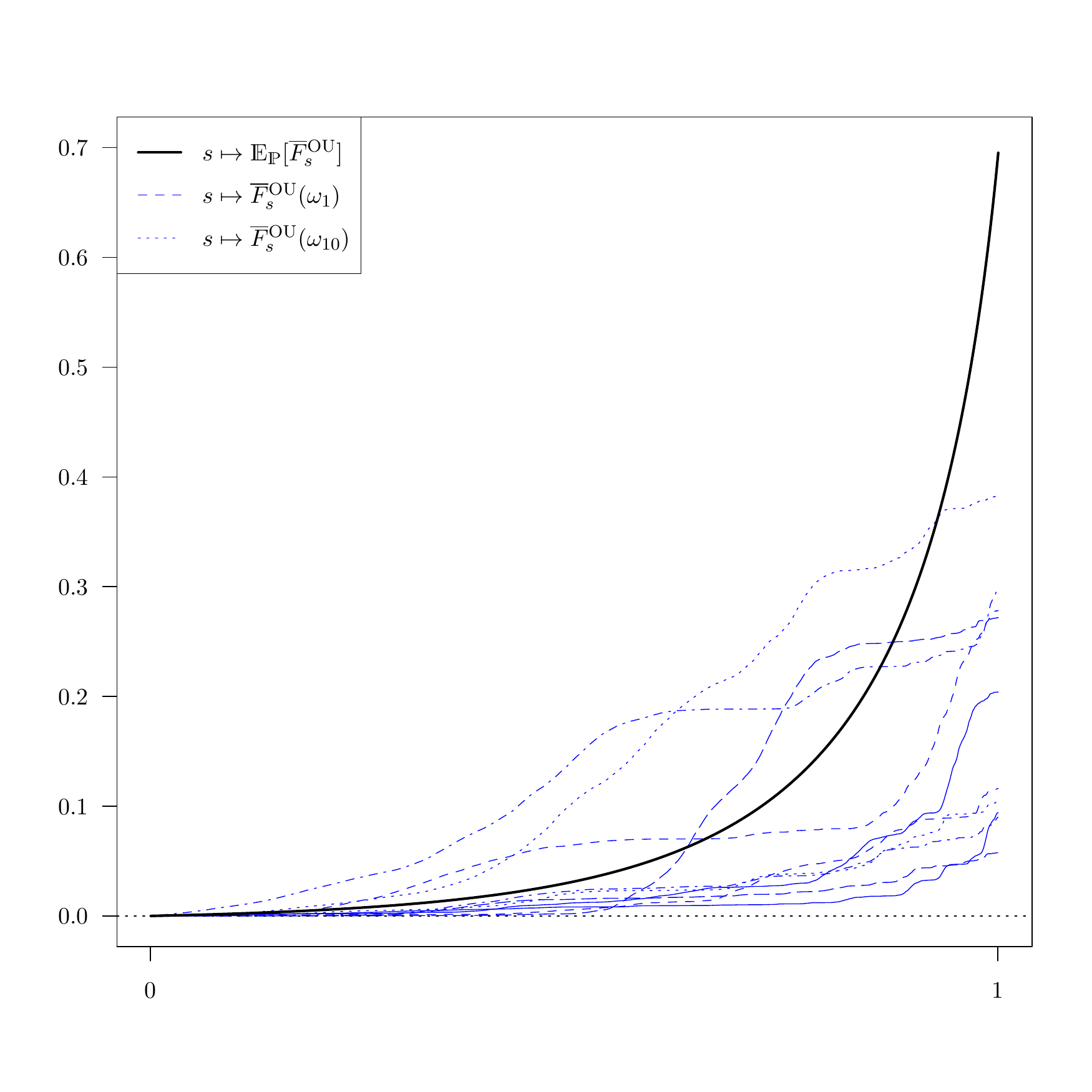}
\caption{Simulations of the cumulative Fisher information process \hyperref[eq: o.u.f.]{(\ref*{eq: o.u.f.})} \\ for the Ornstein--Uhlenbeck diffusion \hyperref[eq: o.u.]{(\ref*{eq: o.u.})}.}
\label{fig: OU.sim}
\end{figure}
\end{example}

\begin{example} We set again $n = 1$ and now consider the case of no confinement potential $V \equiv 0$, quadratic interaction potential $W(x) = \frac{x^{2}}{2}$, and a centered Gaussian initial position $X_{0}$ with variance $\sigma_{0}^{2} > 0$, which is independent of the Brownian motion $(B_{t})_{t \geqslant 0}$ . In this case, for any $t \geqslant 0$, the drift term of the SDE of \hyperref[eq: intro.1.1]{(\ref*{eq: intro.1.1})} is
\begin{equation}
- \nabla (W \ast \mathrm{P}_{\kern-0.1em t})(X_{t})  
= - \int_{\mathds{R}^{n}} \nabla W(X_{t}-y)p_{t}(y) \, \mathrm{d}y
= - \big(X_{t} - \mathds{E}[X_{t}]\big).
\end{equation}
In particular, the drift term depends on the distribution $\mathrm{P}_{\kern-0.1em t}$ only through its mean. Substituting it into \hyperref[eq: intro.1.1]{(\ref*{eq: intro.1.1})}, 
this SDE reduces to
\begin{equation} \label{eq: nonlinear SDE}
\mathrm{d} X_{t} = - \big(X_{t}-\mathds{E}[X_{t}]\big) \, \mathrm{d}t + \sqrt{2} \, \mathrm{d}B_{t} \, , \qquad 0 \leqslant t \leqslant T.
\end{equation}
This type of nonlinear, self-interacting SDE has been studied since \cite{BRTV98}, where it was shown that its solution is also given by the Ornstein--Uhlenbeck process of \hyperref[eq: sol.ou.pr]{(\ref*{eq: sol.ou.pr})}. Therefore, similar computations as in \hyperref[exa: o.u.]{Example \ref*{exa: o.u.}} show that in this setting the cumulative Fisher information process is given by  
\begin{equation}  \label{eq: fis_inf_cw_1}
\mybar{F}_{s}^{\mathrm{NL}} 
= \int_{0}^{s} \Big( \big( \tfrac{1}{\overline{\sigma}_{u}^{4}}  + \tfrac{1}{2}  - \tfrac{1}{\overline{\sigma}_{u}^{2}}  \big) \mybar{X}_{u}^{2}
+ \mathds{E}_{\tilde{\mathds{P}}}\Big[  \tfrac{1}{2} (\mybar{X}_{u}-\mybar{Y}_{u})  \big( \tfrac{2}{\overline{\sigma}_{u}^{2}}  \mybar{Y}_{u} - \mybar{Y}_{u} \big) \Big] \Big) \, \mathrm{d}u
\end{equation}
for $0 \leqslant s \leqslant T$. Using the fact that $(\mybar{X}_{u})_{\#}(\tilde{\mathds{P}}) = (\mybar{Y}_{u})_{\#}(\tilde{\mathds{P}})= \mathcal{N}(0,\overline{\sigma}_{u}^{2})$, which we have from \hyperref[eq: o.u.p.d.f.]{(\ref*{eq: o.u.p.d.f.})}, we can compute the expectation appearing in \hyperref[eq: fis_inf_cw_1]{(\ref*{eq: fis_inf_cw_1})} and obtain
\begin{equation}  \label{eq: fis_inf_cw_2}
\mybar{F}_{s}^{\mathrm{NL}} 
= \int_{0}^{s} \Big( \big( \tfrac{1}{\overline{\sigma}_{u}^{4}}  + \tfrac{1}{2}  - \tfrac{1}{\overline{\sigma}_{u}^{2}}  \big) \mybar{X}_{u}^{2}
+ \tfrac{\overline{\sigma}_{u}^{2}}{2}  -1  \Big) \, \mathrm{d}u.
\end{equation}
Clearly,  $\mybar{F}_{s}^{\mathrm{OU}} \neq \mybar{F}_{s}^{\mathrm{NL}}$. In particular, the integrand in \hyperref[eq: fis_inf_cw_2]{(\ref*{eq: fis_inf_cw_2})} is non-negative if and only if
$\mybar{X}_{u}^{2} \geqslant \big( \tfrac{1}{\overline{\sigma}_{u}^{4}}  + \tfrac{1}{2}  - \tfrac{1}{\overline{\sigma}_{u}^{2}}  \big)^{-1} \big(1 - \tfrac{\overline{\sigma}_{u}^{2}}{2}\big).$
In other words, as opposed to \hyperref[exa: o.u.]{Example \ref*{exa: o.u.}}, relative entropy only decreases along a trajectory if $\mybar{X}_{u}$ is far from its mean.  However, after taking expectations in \hyperref[eq: o.u.f.]{(\ref*{eq: o.u.f.})} and \hyperref[eq: fis_inf_cw_2]{(\ref*{eq: fis_inf_cw_2})}, we see that the expected rate of relative entropy dissipation in both cases is equal to
\begin{equation} \label{eq: exp.path}
\mathds{E}_{\mathds{P}}\big[\mybar{F}_{s}^{\mathrm{OU}}\big] = \mathds{E}_{\mathds{P}}\big[\mybar{F}_{s}^{\mathrm{NL}} \big] = \int_{0}^{s}  \big(\overline{\sigma}_{u}-\tfrac{1}{\overline{\sigma}_{u}}\big)^{2}  \, \mathrm{d}u
\, , \qquad 0 \leqslant s \leqslant T.
\end{equation}

Now we set again $T = 1$ and $\sigma_{0}^{2} = 0.1$. In the same vein as in \hyperref[fig: OU.sim]{Figure \ref*{fig: OU.sim}}, we plot in \hyperref[fig: NL.sim]{Figure \ref*{fig: NL.sim}} below the paths of ten simulated trajectories $s \mapsto \mybar{F}_{s}^{\mathrm{NL}}(\omega_{i})$, for $i = 1, \ldots, 10$. We observe that some of the red lines describing the paths of these trajectories indeed take negative values. In other words, the cumulative Fisher information process of \hyperref[eq: fis_inf_cw_2]{(\ref*{eq: fis_inf_cw_2})}, and hence its integrand, can both be negative. Finally, the thick black line in \hyperref[fig: NL.sim]{Figure \ref*{fig: NL.sim}} follows the expected path $s \mapsto \mathds{E}_{\mathds{P}}[\mybar{F}_{s}^{\mathrm{NL}}]$ of all possible trajectories. According to \hyperref[eq: exp.path]{(\ref*{eq: exp.path})}, this is the same black line as in \hyperref[fig: OU.sim]{Figure \ref*{fig: OU.sim}}.

\vspace*{-1cm}

\begin{figure}[H] 
\centering
\captionsetup{justification=centering}
\includegraphics[width=0.6\textwidth]{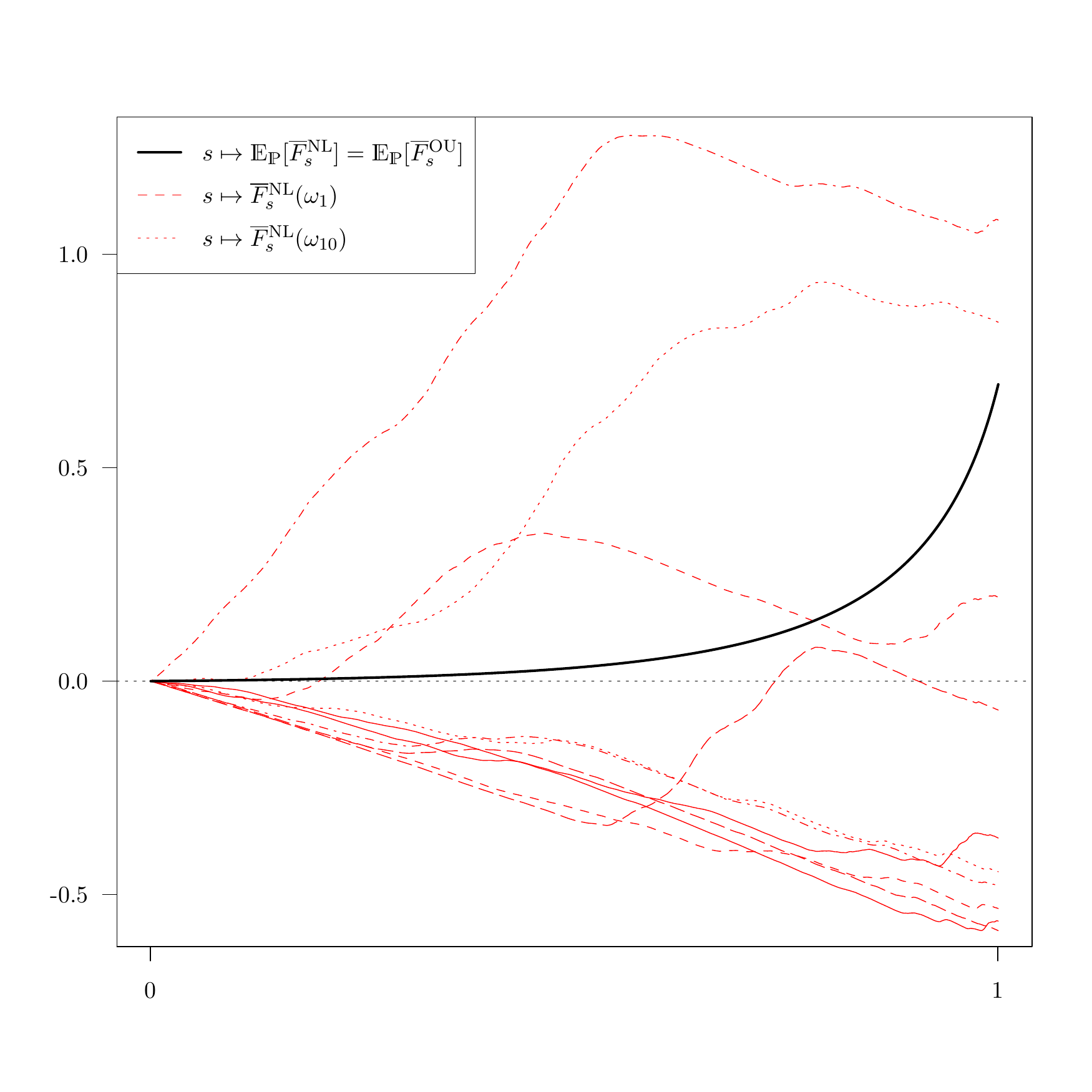}
\caption{Simulations of the cumulative Fisher information process \hyperref[eq: fis_inf_cw_2]{(\ref*{eq: fis_inf_cw_2})} \\ for the nonlinear, self-interacting diffusion \hyperref[{eq: nonlinear SDE}]{(\ref*{eq: nonlinear SDE})}.}
\label{fig: NL.sim}
\end{figure}
\end{example}


\subsubsection{Consequences of \texorpdfstring{\hyperref[thm: trac]{Theorem \ref*{thm: trac}}}{Theorem 3.1}}


We now return to the general statement of \hyperref[thm: trac]{Theorem \ref*{thm: trac}} and deduce some direct consequences. By averaging the trajectorial result of \hyperref[thm: trac]{Theorem \ref*{thm: trac}} according to the path measure $\mathds{P}$, we derive the well-known \emph{relative entropy identity} \hyperref[eq: unper.rel.ent.id.]{(\ref*{eq: unper.rel.ent.id.})} and the \emph{dissipation of relative entropy} \myeqref{eq: unper.rel.ent.dis.} below. A sketch of proof for the latter result was first given in \cite[Proposition 2.1]{CMV03}.

\begin{corollary} \label{cor: ent.diss.id} Suppose \textnormal{\myassumption{assumptions}} hold. For all $t, t_{0} \in [0, T]$, we have the relative entropy identity
\begin{equation} \label{eq: unper.rel.ent.id.} 
H(\mathrm{P}_{\kern-0.1em t} \, \vert \, \mathrm{Q}_{t}) - H(\mathrm{P}_{\kern-0.1em t_{0}} \, \vert \, \mathrm{Q}_{t_{0}}) 
= - \int_{t_{0}}^{t}  I(\mathrm{P}_{\kern-0.1em u} \, \vert \, \mathrm{Q}_{u}^{\uparrow}) \, \textnormal{d}u.
\end{equation}
In particular, the relative entropy function $t \mapsto H(\mathrm{P}_{\kern-0.1em t} \, \vert \, \mathrm{Q}_{t})$ is monotonically decreasing. Furthermore, for Lebesgue-a.e.\ $t \in [0,T]$, the relative Fisher information $I(\mathrm{P}_{\kern-0.1em t} \, \vert \, \mathrm{Q}_{t}^{\uparrow})$ is finite, and the rate of relative entropy dissipation is given by
\begin{equation} \label{eq: unper.rel.ent.dis.} 
\tfrac{\textnormal{d}}{\textnormal{d}t} \, H(\mathrm{P}_{\kern-0.1em t} \, \vert \, \mathrm{Q}_{t}) 
= -  I(\mathrm{P}_{\kern-0.1em t} \, \vert \, \mathrm{Q}_{t}^{\uparrow}).
\end{equation}
\begin{proof} The identity \hyperref[eq: unper.rel.ent.id.]{(\ref*{eq: unper.rel.ent.id.})} follows by taking expectations with respect to the probability measure $\mathds{P}$ in \hyperref[eq: d.m.d]{(\ref*{eq: d.m.d})}, recalling \hyperref[eq: intro.1.8]{(\ref*{eq: intro.1.8})}, using \hyperref[eq: fish.inf.exp]{(\ref*{eq: fish.inf.exp})}, and invoking the fact that the $\mathds{P}$-expectation of the martingale \hyperref[eq: martingale]{(\ref*{eq: martingale})} vanishes. Finally, applying the Lebesgue differentiation theorem to the monotone function $t \mapsto H(\mathrm{P}_{\kern-0.1em t} \, \vert \, \mathrm{Q}_{t})$ gives \myeqref{eq: unper.rel.ent.dis.}. 
\end{proof}
\end{corollary}

\begin{remark} The relation \hyperref[eq: unper.rel.ent.dis.]{(\ref*{eq: unper.rel.ent.dis.})} describes the temporal dissipation of relative entropy at the \emph{ensemble} level. It asserts that the rate of decay of the relative entropy $t \mapsto H(\mathrm{P}_{\kern-0.1em t} \, \vert \, \mathrm{Q}_{t})$ is given by the relative Fisher information $I(\mathrm{P}_{\kern-0.1em t} \, \vert \, \mathrm{Q}_{t}^{\uparrow})$. 
\end{remark}

Finally, let us place ourselves again on the filtered probability space $(\Omega,\mathds{G},\mathds{P})$ as in \hyperref[thm: trac]{Theorem \ref*{thm: trac}} and emphasize that this trajectorial result is valid along almost every trajectory $s \mapsto \mybar{X}_{s}(\omega)$ of the underlying McKean--Vlasov process. As a consequence, by taking conditional expectations, we can generalize \myeqref{eq: unper.rel.ent.dis.} and deduce the following trajectorial rate of relative entropy dissipation. 

\begin{corollary} \label{cor: traj.de.Bruijn} Suppose \textnormal{\myassumption{assumptions}} hold and $\int_{0}^{T} \mathds{E}_{\mathds{P}} [ \vert \mybar{I}_{u} \vert ] \, \mathrm{d}u  < \infty$. For $\mathds{P}$-a.e.\ $\omega \in \Omega$ there exists a Lebesgue null set $N_{\omega} \subseteq [0,T]$ such that for any $s_{0} \in [0,T] \setminus N_{\omega}$ we have
\begin{equation} \label{eq: traj.de.Bruijn.a} 
\lim_{s \downarrow s_{0}} \frac{\mathds{E}_{\mathds{P}} \big[\log \mybar{\ell}_{s}(\mybar{X}_{s},\mybar{\mathrm{P}}_{\kern-0.1em s})  \, \vert \, \mathcal{G}_{s_{0}}\big](\omega) - \log \mybar{\ell}_{s_{0}}\big(\mybar{X}_{s_{0}}(\omega),\mybar{\mathrm{P}}_{\kern-0.1em s_{0}}\big)}{s-s_{0}} 
=  \mybar{I}_{s_{0}}(\omega).
\end{equation}

\begin{remark} Recalling \hyperref[eq: fish.inf.exp]{(\ref*{eq: fish.inf.exp})}, we observe that $\mathds{E}_{\mathds{P}}[\mybar{I}_{s_{0}}] = I(\mybar{\mathrm{P}}_{\kern-0.1em s_{0}} \, \vert \, \mybar{\mathrm{Q}}_{s_{0}}^{\uparrow})$. Therefore the limiting assertion \hyperref[eq: traj.de.Bruijn.a]{(\ref*{eq: traj.de.Bruijn.a})} can indeed be viewed as a trajectorial version of the deterministic relative entropy dissipation identity \hyperref[eq: unper.rel.ent.dis.]{(\ref*{eq: unper.rel.ent.dis.})}. 
\end{remark}

\begin{proof} We let $0 \leqslant s_{0} \leqslant s \leqslant T$. By \hyperref[eq: d.m.d]{(\ref*{eq: d.m.d})}, \hyperref[eq: c.f.i.]{(\ref*{eq: c.f.i.})} and Fubini's theorem, we have for $\mathds{P}$-a.e.\ $\omega \in \Omega$ 
\begin{equation} \label{eq: ljc.a}
\mathds{E}_{\mathds{P}} \big[\log \mybar{\ell}_{s}(\mybar{X}_{s},\mybar{\mathrm{P}}_{\kern-0.1em s})  \, \vert \, \mathcal{G}_{s_{0}}\big](\omega) - \log \mybar{\ell}_{s_{0}}\big(\mybar{X}_{s_{0}}(\omega),\mybar{\mathrm{P}}_{\kern-0.1em s_{0}}\big) 
= \int_{s_{0}}^{s} \mathds{E}_{\mathds{P}} \big[ \mybar{I}_{u} \, \vert \, \mathcal{G}_{s_{0}} \big](\omega) \,  \mathrm{d}u.
\end{equation}
Furthermore, Jensen's inequality gives  
\begin{equation}
\mathds{E}_{\mathds{P}} \bigg[ \int_{s_{0}}^{s} \big\vert \mathds{E}_{\mathds{P}} \big[ \mybar{I}_{u} \, \vert \, \mathcal{G}_{s_{0}} \big] \big\vert \,  \mathrm{d}u \bigg] 
\leqslant \int_{s_{0}}^{s} \mathds{E}_{\mathds{P}} \big[ \vert \mybar{I}_{u} \vert \big] \, \mathrm{d}u 
< \infty, 
\end{equation}
which implies 
\begin{equation}
\int_{s_{0}}^{s} \big\vert \mathds{E}_{\mathds{P}} \big[ \mybar{I}_{u} \, \vert \, \mathcal{G}_{s_{0}} \big] (\omega) \big\vert \,  \mathrm{d}u < \infty
\qquad \textnormal{ for } \mathds{P}\textnormal{-a.e. } \omega \in \Omega.
\end{equation}
By the Lebesgue differentiation theorem, for every such $\omega$ there exists a Lebesgue null set $N_{\omega} \subseteq [0,T]$ so that the limiting assertion 
\begin{equation} \label{eq: ljc.b}
\lim_{s \downarrow s_{0}} \frac{\int_{s_{0}}^{s} \mathds{E}_{\mathds{P}} \big[ \mybar{I}_{u} \, \vert \, \mathcal{G}_{s_{0}} \big](\omega) \,  \mathrm{d}u}{s-s_{0}} = \mathds{E}_{\mathds{P}} \big[ \mybar{I}_{s_{0}} \, \vert \, \mathcal{G}_{s_{0}} \big](\omega)  =  \mybar{I}_{s_{0}}(\omega)
\end{equation}
holds for every $s_{0} \in [0,T] \setminus N_{\omega}$. Finally, combining \hyperref[eq: ljc.a]{(\ref*{eq: ljc.a})} and \hyperref[eq: ljc.b]{(\ref*{eq: ljc.b})} proves \hyperref[eq: traj.de.Bruijn.a]{(\ref*{eq: traj.de.Bruijn.a})}.
\end{proof}
\end{corollary}


\subsection{Gradient flow structure of the granular media equation} \label{subsec: gradient flow}


In this subsection we apply the trajectorial approach of \mysubsecref{subsec: traj diss} in order to formulate the gradient flow property of the granular media equation \myeqref{eq: intro.1.2}. To this end, we consider a function $\boldsymbol{\beta} \colon \mathds{R}^{n} \rightarrow \mathds{R}$, which will be treated as a perturbation potential. We denote by $V^{\boldsymbol{\beta}} \coloneqq V + \boldsymbol{\beta}$ the perturbed confinement potential and invoke the following regularity assumptions.

\begin{assumptions} \label{assumptions.b} The function $\bbeta \colon \mathds{R}^{n} \rightarrow \mathds{R}$ is smooth and compactly supported, and we require that \myassumption{assumptions} are still satisfied if we replace $V$ by $V^{\bbeta}$.
\end{assumptions}

Note that \myassumption{assumptions.b} allow us to apply \mylemref{lem: well-posedness} to the ``perturbed'' McKean--Vlasov SDE 
\begin{equation} \label{eq: perturbed SDE}
\mathrm{d}X_{t} = - \big( \nabla V^{\boldsymbol{\beta}}(X_{t}) + \nabla (W \ast \mathrm{P}_{\kern-0.1em t}^{\boldsymbol{\beta}})(X_{t}) \big) \, \mathrm{d}t + \sqrt{2} \, \mathrm{d}B_{t}^{\boldsymbol{\beta}} 
\, , \qquad t_{0} \leqslant t \leqslant T
\end{equation}
starting at time $t_{0} \geqslant 0$, with $X_{t_{0}}$ having initial distribution $\mathrm{P}_{\kern-0.1em t_{0}}^{\boldsymbol{\beta}} = \mathrm{P}_{\kern-0.1em t_{0}}$. Therefore, by analogy with \mysubsecref{sec: notations}, we can construct a probability measure $\mathds{P}^{\boldsymbol{\beta}}$ on $\Omega \coloneqq C([t_{0},T];\mathds{R}^{n})$, under which the canonical process $(X_{t})_{t_{0} \leqslant t \leqslant T}$ satisfies the SDE \myeqref{eq: perturbed SDE}, with $(B_{t}^{\boldsymbol{\beta}})_{t_{0} \leqslant t \leqslant T}$ being a $\P^{\bbeta}$-Brownian motion.

For each time $t \in [t_{0},T]$, we denote by $\mathrm{P}_{\kern-0.1em t}^{\boldsymbol{\beta}} \coloneqq \P^{\bbeta}\circ X_t^{-1}$ the probability distribution and by $p_{t}^{\boldsymbol{\beta}}$ the probability density function of $X_{t}$ under $\mathds{P}^{\boldsymbol{\beta}}$. The ``perturbed'' curve of density functions $(p_{t}^{\boldsymbol{\beta}})_{t_{0} \leqslant t \leqslant T}$ then satisfies the \textit{perturbed granular media equation}
\begin{equation} \label{eq: perturbed kinetic FPK}
\begin{cases}
\partial_{t} p_{t}^{\boldsymbol{\beta}}(x) = \operatorname{div} \Big( \nabla p_{t}^{\boldsymbol{\beta}}(x) + p_{t}^{\boldsymbol{\beta}}(x)  \nabla V^{\boldsymbol{\beta}}(x) + p_{t}^{\boldsymbol{\beta}}(x)  \nabla(W \ast p_{t}^{\boldsymbol{\beta}})(x)\Big) \, , & (t,x) \in (t_{0},T) \times \mathds{R}^{n},    \\[7pt]
p_{t_{0}}^{\boldsymbol{\beta}}(x) = p_{t_{0}}(x) \, , & x \in \mathds{R}^{n}.
\end{cases}
\end{equation}

By analogy with \hyperref[eq: pot]{(\ref*{eq: pot})}, we define the perturbed potentials 
\begin{equation} \label{eq: perturbed shifted potential}
\Psi^{\boldsymbol{\beta}}(x,\mu) \coloneqq V^{\boldsymbol{\beta}}(x) + \tfrac{1}{2}(W\ast\mu)(x)
\, , \qquad
\Psi^{\boldsymbol{\beta}\uparrow}(x,\mu) \coloneqq V^{\boldsymbol{\beta}}(x) + (W\ast\mu)(x)
\, , \qquad
\Psi^{\boldsymbol{\beta}\downarrow} \coloneqq V^{\boldsymbol{\beta}}
\end{equation}
for $(x,\mu) \in \mathds{R}^{n} \times \mathscr{P}_{2}(\mathds{R}^{n})$. In parallel to \hyperref[eq: intro.1.7]{(\ref*{eq: intro.1.7})}, we introduce the \textit{perturbed likelihood ratio functions} 
\begin{equation} \label{eq: perturbed  likeilhood rato, perturbed  stationary distribution}
\ell_{t}^{\boldsymbol{\beta}}(x,\mu) \coloneqq \frac{p_{t}^{\boldsymbol{\beta}}(x)}{q(x,\mu)} 
\ , \qquad \qquad
\ell_{t}^{\boldsymbol{\beta}\uparrow}(x,\mu) \coloneqq \frac{p_{t}^{\boldsymbol{\beta}}(x)}{q^{\uparrow}(x,\mu)} 
\ , \qquad \qquad
\ell_{t}^{\boldsymbol{\beta}\downarrow}(x) \coloneqq \frac{p_{t}^{\boldsymbol{\beta}}(x)}{q^{\downarrow}(x)}
\end{equation} 
for $t \in [t_{0},T]$. Finally, we define the $\sigma$-finite measures
\begin{equation} 
\mathrm{Q}_{t}^{\boldsymbol{\beta}}(A) \coloneqq \int_{A} q(x,\mathrm{P}_{\kern-0.1em t}^{\boldsymbol{\beta}}) \, \mathrm{d}x 
\, , \qquad 
\mathrm{Q}_{t}^{\boldsymbol{\beta}\uparrow}(A) \coloneqq \int_{A} q^{\uparrow}(x,\mathrm{P}_{\kern-0.1em t}^{\boldsymbol{\beta}}) \, \mathrm{d}x \, , \qquad A \in \mathscr{B}(\mathds{R}^{n}).
\end{equation}
They are the perturbed versions of the measures $\mathrm{Q}_{t}$ and $\mathrm{Q}_{t}^{\uparrow}$ defined in \hyperref[eq: intro.1.6]{(\ref*{eq: intro.1.6})}. The relative entropy of $\mathrm{P}_{\kern-0.1em t}^{\boldsymbol{\beta}}$ with respect to $\mathrm{Q}_{t}^{\boldsymbol{\beta}}$ is then given by
\begin{equation} \label{eq: per.re.ent}
H\big(\mathrm{P}_{\kern-0.1em t}^{\boldsymbol{\beta}} \, \big\vert \, \mathrm{Q}_{t}^{\boldsymbol{\beta}}\big) = 
\mathds{E}_{\mathds{P}^{\boldsymbol{\beta}}}\big[\log\ell_{t}^{\boldsymbol{\beta}}(X_{t},\mathrm{P}_{\kern-0.1em t}^{\boldsymbol{\beta}})\big] = 
\mathscr{F}(\mathrm{P}_{\kern-0.1em t}^{\boldsymbol{\beta}})
\end{equation}
and the relative Fisher information of $\mathrm{P}_{\kern-0.1em t}^{\boldsymbol{\beta}}$ with respect to $\mathrm{Q}_{t}^{\boldsymbol{\beta}\uparrow}$ equals
\begin{equation} \label{eq: pe.fi.inf.id}
I\big(\mathrm{P}_{\kern-0.1em t}^{\boldsymbol{\beta}} \, \big\vert \, \mathrm{Q}_{t}^{\boldsymbol{\beta}\uparrow}\big) = \mathds{E}_{\mathds{P}^{\boldsymbol{\beta}}}\Big[\big\vert\nabla\log\ell_{t}^{\boldsymbol{\beta}\uparrow}(X_{t},\mathrm{P}_{\kern-0.1em t}^{\boldsymbol{\beta}})\big\vert^{2}\Big] =
\Ds(\mathrm{P}_{\kern-0.1em t}^{\boldsymbol{\beta}}).
\end{equation}

The following trajectorial result, \hyperref[thm: p.trac]{Theorem \ref*{thm: p.trac}} below, provides the semimartingale decomposition of the \textit{perturbed relative entropy process}
\begin{equation} \label{eq: p.r.e.p.}
\log \ell_{t}^{\boldsymbol{\beta}} (X_{t},\mathrm{P}_{\kern-0.1em t}^{\boldsymbol{\beta}}) = \log p_{t}^{\boldsymbol{\beta}}(X_{t}) + V(X_{t}) + \tfrac{1}{2}(W \ast \mathrm{P}_{\kern-0.1em t}^{\boldsymbol{\beta}})(X_{t}) 
\, , \qquad t_{0} \leqslant t \leqslant T. 
\end{equation}
In line with its unperturbed counterpart, \hyperref[thm: trac]{Theorem \ref*{thm: trac}}, we shall formulate this result in the reverse direction of time. We first introduce the perturbed analogues of \hyperref[eq: f.i.p.]{(\ref*{eq: f.i.p.})} and \hyperref[eq: c.f.i.]{(\ref*{eq: c.f.i.})}: the \emph{time-reversed perturbed Fisher information process} is defined as
\begingroup
\addtolength{\jot}{0.7em} 
\begin{align}
\mybar{I}_{s}^{\boldsymbol{\beta}} &\coloneqq   \Big( \big\vert \nabla\log\mybar{\ell}_{s}^{\boldsymbol{\beta}\downarrow} \big\vert^{2} + \tfrac{1}{2}  \big\vert \nabla( W\ast\mybar{\mathrm{P}}_{\kern-0.1em s}^{\boldsymbol{\beta}}) \big\vert^{2} + \Big\langle  \tfrac{1}{2} \nabla(W\ast\mybar{\mathrm{P}}_{\kern-0.1em s}^{\boldsymbol{\beta}}) \,  , \, 2\nabla\log\mybar{\ell}_{s}^{\boldsymbol{\beta}\downarrow} + \nabla V^{\boldsymbol{\beta}}  \Big\rangle \Big)(\mybar{X}_{s}) \,  \label{eq: p.c.f.i.a} \\
& \qquad  -  \mathds{E}_{\tilde{\mathds{P}}^{\boldsymbol{\beta}}}\Big[\Big\langle  \tfrac{1}{2} \nabla W(\mybar{X}_{s}-\mybar{Y}_{s})  \, , \Big(2\nabla\log\mybar{\ell}_{s}^{\boldsymbol{\beta}\downarrow} - \nabla V +\nabla ( W\ast\mybar{\mathrm{P}}_{\kern-0.1em s}^{\boldsymbol{\beta}}) + \nabla \boldsymbol{\beta} \Big)(\mybar{Y}_{s})\Big\rangle \Big] \,  \label{eq: p.c.f.i.b} \\
&\qquad + \Big(  \langle \nabla V \, , \nabla \boldsymbol{\beta} \rangle -\Delta \boldsymbol{\beta} \Big)(\mybar{X}_{s}) \label{eq: p.c.f.i.c} 
\end{align}
\endgroup
for all $0 \leqslant s \leqslant T-t_{0}$, where $(\mybar{Y}_{s})_{0 \leqslant s \leqslant T-t_{0}}$ is a copy of the process $(\mybar{X}_{s})_{0 \leqslant s \leqslant T-t_{0}}$ on a copy $(\tilde{\Omega},\tilde{\G},\tilde{\mathds{P}}^{\boldsymbol{\beta}})$ of the original probability space $(\Omega,\G,\mathds{P}^{\boldsymbol{\beta}})$; the \emph{perturbed cumulative Fisher information process} is defined as 
\begin{equation} \label{eq: pcfip}
\mybar{F}_{s}^{\boldsymbol{\beta}} \coloneqq \int_{0}^{s} \mybar{I}_{u}^{\boldsymbol{\beta}} \, \mathrm{d}u  \, , \qquad 0 \leqslant s \leqslant T-t_{0}.
\end{equation}

\begin{theorem} \label{thm: p.trac} Suppose \textnormal{\myassumption{assumptions.b}} hold. On $(\Omega,\mathds{G},\P^{\bbeta})$, the time-reversed perturbed relative entropy process 
\begin{equation} \label{eq: trprep}
\log \mybar{\ell}_{s}^{\boldsymbol{\beta}}(\mybar{X}_{s},\mybar{\mathrm{P}}_{\kern-0.1em s}^{\boldsymbol{\beta}}) = \log\ell_{T-s}^{\boldsymbol{\beta}}(X_{T-s},\mathrm{P}_{\kern-0.1em T-s}^{\boldsymbol{\beta}}) 
\, , \qquad  0 \leqslant s \leqslant T-t_{0}
\end{equation}
admits the semimartingale decomposition 
\begin{equation} \label{eq: p.dmd}
\log \mybar{\ell}_{s}^{\boldsymbol{\beta}}(\mybar{X}_{s},\mybar{\mathrm{P}}_{\kern-0.1em s}^{\boldsymbol{\beta}}) - \log \mybar{\ell}_{0}^{\boldsymbol{\beta}}(\mybar{X}_{0},\mybar{\mathrm{P}}_{\kern-0.1em 0}^{\boldsymbol{\beta}}) = \mybar{M}_{s}^{\boldsymbol{\beta}} + \mybar{F}_{s}^{\boldsymbol{\beta}}.
\end{equation}
Here $(\mybar{M}_{s}^{\bbeta})_{0 \leqslant s \leqslant T-t_{0}}$ is the $L^{2}(\mathds{P}^{\bbeta})$-bounded martingale 
\begin{equation} \label{eq: p.mart}
\mybar{M}_{s}^{\boldsymbol{\beta}} \coloneqq \int_{0}^{s} \Big\langle \nabla \log \mybar{\ell}_{u}^{\boldsymbol{\beta}}(\mybar{X}_{u},\mybar{\mathrm{P}}_{\kern-0.1em u}^{\boldsymbol{\beta}}) \, , \, \sqrt{2} \, \mathrm{d}\mybar{B}_{u}^{\boldsymbol{\beta}} \Big\rangle,
\end{equation}
with $(\mybar{B}_{s}^{\bbeta})_{0 \leqslant s \leqslant T-t_{0}}$ a $\P^{\bbeta}$-Brownian motion of the backward filtration $\G$, and the compensator \textnormal{\hyperref[eq: pcfip]{(\ref*{eq: pcfip})}} satisfies
\begin{equation} \label{eq: fish.inf.exp.pert}
\mathds{E}_{\mathds{P}^{\boldsymbol{\beta}}}\big[\mybar{F}_{s}^{\boldsymbol{\beta}}\big] 
= \int_{0}^{s} \bigg( I\big(\mybar{\mathrm{P}}_{\kern-0.1em u}^{\boldsymbol{\beta}} \, \big\vert \, \mybar{\mathrm{Q}}_{u}^{\boldsymbol{\beta}\uparrow}\big) + \mathds{E}_{\mathds{P}^{\boldsymbol{\beta}}} \Big[ \Big(  \Big\langle \nabla V + \nabla( W\ast\mybar{\mathrm{P}}_{\kern-0.1em u}^{\boldsymbol{\beta}}) \, , \nabla \boldsymbol{\beta} \Big\rangle -\Delta \boldsymbol{\beta} \Big)(\mybar{X}_{u})\Big] \bigg) \, \mathrm{d}u < \infty.
\end{equation}
\end{theorem}

With the dynamics of the time-reversed perturbed relative entropy process at hand, we repeat the same procedure which was carried out for the unperturbed case. Taking expectations with respect to the probability measure $\mathds{P}^{\boldsymbol{\beta}}$, we arrive at the  perturbed relative entropy identity \myeqref{eq: per.rel.ent.id.}, and applying the Lebesgue differentiation theorem gives the perturbed relative entropy production identity \hyperref[eq: per.rel.ent.dis.]{(\ref*{eq: per.rel.ent.dis.})}.

\begin{corollary} \label{cor: per.rel.ent.id.b} Suppose \textnormal{\myassumption{assumptions.b}} hold. For all $0 \leqslant t_{0} \leqslant t \leqslant T$, we have the perturbed relative entropy identity
\begin{equation} \label{eq: per.rel.ent.id.} 
H\big(\mathrm{P}_{\kern-0.1em t}^{\boldsymbol{\beta}} \, \big\vert \, \mathrm{Q}_{t}^{\boldsymbol{\beta}}\big) - H\big(\mathrm{P}_{\kern-0.1em t_{0}}^{\boldsymbol{\beta}} \, \big\vert \, \mathrm{Q}_{t_{0}}^{\boldsymbol{\beta}}\big) 
= - \int_{t_{0}}^{t} \bigg( I\big(\mathrm{P}_{\kern-0.1em u}^{\boldsymbol{\beta}} \, \big\vert \, \mathrm{Q}_{u}^{\boldsymbol{\beta}\uparrow}\big) + \mathds{E}_{\mathds{P}^{\boldsymbol{\beta}}} \Big[ \Big(  \Big\langle \nabla V + \nabla( W\ast\mathrm{P}_{\kern-0.1em u}^{\boldsymbol{\beta}}) \, , \nabla \boldsymbol{\beta} \Big\rangle -\Delta \boldsymbol{\beta} \Big)(X_{u})\Big] \bigg) \, \textnormal{d}u.
\end{equation}
For Lebesgue-a.e.\ $t_{0} \in [0, T]$, the perturbed rate of relative entropy dissipation is given by
\begin{equation} \label{eq: per.rel.ent.dis.} 
\frac{\mathrm{d}}{\mathrm{d}t} \Big\vert_{t=t_{0}}^{+} \,  H\big(\mathrm{P}_{\kern-0.1em t}^{\boldsymbol{\beta}} \, \big\vert \, \mathrm{Q}_{t}^{\boldsymbol{\beta}}\big) 
= -  \bigg( I\big(\mathrm{P}_{\kern-0.1em t_{0}} \, \big\vert \, \mathrm{Q}_{t_{0}}\big) + \mathds{E}_{\mathds{P}} \Big[ \Big(  \Big\langle \nabla V + \nabla( W\ast\mathrm{P}_{\kern-0.1em t_{0}}) \, , \nabla \boldsymbol{\beta} \Big\rangle - \Delta \boldsymbol{\beta} \Big)(X_{t_{0}})\Big] \bigg).
\end{equation}
\end{corollary}

Similarly, we have the following trajectorial rate of relative entropy dissipation for the perturbed diffusion.

\begin{corollary} \label{cor: p.traj.de.Bruijn} Suppose \textnormal{\myassumption{assumptions.b}} hold and $\int_{0}^{T-t_{0}} \mathds{E}_{\mathds{P}^{\boldsymbol{\beta}}} [ \vert \mybar{I}_{u}^{\boldsymbol{\beta}} \vert ] \, \mathrm{d}u  < \infty$. For $\mathds{P}^{\boldsymbol{\beta}}$-a.e.\ $\omega \in \Omega$ there exists a Lebesgue null set $N_{\omega}^{\boldsymbol{\beta}} \subseteq [0,T-t_{0}]$ such that for any $s_{0} \in [0,T-t_{0}] \setminus N_{\omega}^{\boldsymbol{\beta}}$ we have
\begin{equation} \label{eq: p.traj.de.Bruijn.a}
\lim_{s \downarrow s_{0}} \frac{\mathds{E}_{\mathds{P}^{\boldsymbol{\beta}}} \big[\log \mybar{\ell}_{s}^{\boldsymbol{\beta}}(\mybar{X}_{s},\mybar{\mathrm{P}}_{\kern-0.1em s}^{\boldsymbol{\beta}})  \, \vert \, \mathcal{G}_{s_{0}}\big](\omega) - \log \mybar{\ell}_{s_{0}}^{\boldsymbol{\beta}}\big(\mybar{X}_{s_{0}}(\omega),\mybar{\mathrm{P}}_{\kern-0.1em s_{0}}^{\boldsymbol{\beta}}\big)}{s-s_{0}}  =  \mybar{I}_{s_{0}}^{\boldsymbol{\beta}}(\omega).
\end{equation}
\begin{proof} The proof proceeds almost verbatim as the proof of \hyperref[cor: traj.de.Bruijn]{Corollary \ref*{cor: traj.de.Bruijn}}. The only difference is that we now use the semimartingale decomposition \hyperref[eq: p.dmd]{(\ref*{eq: p.dmd})} and the $\mathds{P}^{\boldsymbol{\beta}}$-martingale property of the process \hyperref[eq: p.mart]{(\ref*{eq: p.mart})} in \hyperref[thm: p.trac]{Theorem \ref*{thm: p.trac}}. 
\end{proof} 
\end{corollary}

We now turn to the computation of the rate of change of the Wasserstein distance along the curve of probability distributions $(\Pr^{\bbeta}_{\kern-0.1em t})_{t_{0} \leqslant t \leqslant T}$. To this end, we set
\begin{equation} \label{eq: vel.fiel}
v_{t}^{\bbeta}(x) \coloneqq -\big( \nabla \log p_{t}^{\bbeta} + \nabla V^{\bbeta} +  \nabla(W \ast p_{t}^{\bbeta})\big)(x) \, , \qquad (t,x) \in [t_{0},T] \times \mathds{R}^{n},
\end{equation}
so that the perturbed granular media equation \hyperref[eq: perturbed kinetic FPK]{(\ref*{eq: perturbed kinetic FPK})} can be viewed as a \emph{continuity equation}
\begin{equation}
\partial_{t} p_{t}^{\bbeta}(x) + \operatorname{div} \big( v_{t}^{\bbeta}(x) \, p_{t}^{\bbeta}(x) \big) = 0 \, , \qquad (t,x) \in (t_{0},T) \times \mathds{R}^{n},
\end{equation}
with $v_{t}^{\bbeta}(\, \cdot \,)$ as the corresponding \emph{velocity field}. We recall the definition of the tangent space (see Definition 8.4.1 in \cite{AGS08}) 
\begin{equation}
\operatorname{Tan}_{\mu} \mathscr{P}_{2}(\mathds{R}^{n}) \coloneqq \overline{\{\nabla \varphi \colon \varphi \in C_{c}^{\infty}(\mathds{R}^{n};\mathds{R})\}}^{L^{2}(\mu)}
\end{equation}
of $\Ps_{2}(\R^{n})$ at the point $\mu \in \Ps_{2}(\R^{n})$, and impose the following additional assumptions.

\begin{assumptions} \label{assumptions.c} In addition to \myassumption{assumptions.b}, we suppose that 
\begin{equation} \label{eq: tangent.space.assum}
v_{t}(\, \cdot \,) \in \operatorname{Tan}_{\mathrm{P}_{\kern-0.1em t}}  \mathscr{P}_{2}(\mathds{R}^{n}) \qquad \textnormal{ for Lebesgue-a.e.\ } t \in [0,T],
\end{equation}
where $v_{t}(\, \cdot \,)$ is obtained by taking $\bbeta \equiv 0$ and $t_{0} = 0$ in \hyperref[eq: vel.fiel]{(\ref*{eq: vel.fiel})}.
\end{assumptions}

\begin{remark} For example, we know from \cite[Theorem 10.4.13]{AGS08} that the condition \myeqref{eq: tangent.space.assum} is satisfied if, in addition to \myassumption{assumptions.b}, $V$ is uniformly convex, i.e., $\mathrm{Hess}(V) \geqslant \kappa_{V} I_{n}$ for some real constant $\kappa_{V}$, and $W$ is a convex function satisfying the \emph{doubling condition}
\begin{equation}
\exists \, C_{W} > 0 \, \textnormal{ such that } \, \forall \, x, y \in \mathds{R}^{n} \colon \qquad W(x+y) \leqslant C_{W} \big(1 + W(x) + W(y)\big).
\end{equation}
\end{remark}

The proof of the following result is based on the general theory of Wasserstein metric derivatives of absolutely continuous curves in $\mathscr{P}_{\mathrm{ac},2}(\mathds{R}^{n})$; for a thorough discussion, we refer to Chapter 8 in \cite{AGS08}.

\begin{lemma} \label{lem: p.l.b.w.d.} Suppose \textnormal{\myassumption{assumptions.c}} hold. For Lebesgue-a.e.\ $t_{0} \in [0,T]$, the Wasserstein metric derivative of the perturbed curve $(\mathrm{P}_{\kern-0.1em t}^{\bbeta})_{t_{0} \leqslant t \leqslant T}$ is equal to
\begin{equation} \label{eq: local perturbed wasserstein}  
\lim_{t \downarrow t_{0}} \frac{W_{2}\big(\Pr_{\kern-0.1em t}^{\bbeta},\Pr_{\kern-0.1em t_0}^{\bbeta} \big)}{t-t_0}  
= \big\Vert v_{t_{0}}^{\boldsymbol{\beta}}(X_{t_0}) \big\Vert_{L^{2}(\mathds{P})}
= \big\Vert \nabla \log \ell_{t_0}^{\uparrow}(X_{t_0},\mathrm{P}_{\kern-0.1em {t_0}}) + \nabla \boldsymbol{\beta}(X_{t_0}) \big\Vert_{L^{2}(\mathds{P})}.
\end{equation}

\begin{proof} Without loss of generality we can set $\boldsymbol{\beta} \equiv 0$. Note that from \hyperref[eq: fish.inf.exp]{(\ref*{eq: fish.inf.exp})} we have
\begin{equation} 
\int_{0}^{T} \int_{\mathds{R}^{n}} \vert v_{t}(x) \vert^{2} \, \mathrm{d}p_{t}(x) \, \mathrm{d}t 
= \mathds{E}_{\mathds{P}}\bigg[ \int_{0}^{T} \big\vert\nabla\log\ell_{t}^{\uparrow}(X_{t},\mathrm{P}_{\kern-0.1em t})\big\vert^{2} \, \mathrm{d}t \bigg] 
< \infty,
\end{equation}
which implies that $v_{t}(\, \cdot \, ) \in L^{2}(\mathrm{P}_{\kern-0.1em t})$ for Lebesgue-a.e.\ $t \in [0,T]$. Therefore we can apply Theorem 8.3.1 and Proposition 8.4.5 of \cite{AGS08} to the absolutely continuous curve $(\mathrm{P}_{\kern-0.1em t})_{0 \leqslant t \leqslant T}$, which yields \hyperref[eq: local perturbed wasserstein]{(\ref*{eq: local perturbed wasserstein})}.
\end{proof}
\end{lemma}

We now have all the ingredients to formulate the gradient flow property of the granular media equation. The Wasserstein metric slope of the free energy functional $\mathscr{F}$ along the McKean--Vlasov curve $(\mathrm{P}_{\kern-0.1em t})_{t_{0} \leqslant t \leqslant T}$ is defined as 
\begin{equation}
\big\vert \partial \mathscr{F} \big\vert_{W_{2}}(\mathrm{P}_{\kern-0.1em t_{0}}) 
\coloneqq \lim_{t \downarrow t_{0}} \frac{H(\mathrm{P}_{\kern-0.1em t} \, \vert \, \mathrm{Q}_{t}) - H(\mathrm{P}_{\kern-0.1em t_{0}} \, \vert \, \mathrm{Q}_{t_{0}})}{W_{2}(\mathrm{P}_{\kern-0.1em t},\mathrm{P}_{\kern-0.1em t_{0}})}.
\end{equation}
In order to show that this is the slope of steepest descent, we will compare it with the slope 
\begin{equation} \label{eq: p.wass.s}
\big\vert \partial \mathscr{F}^{\boldsymbol{\beta}} \big\vert_{W_{2}}\big(\mathrm{P}_{\kern-0.1em t_{0}}^{\boldsymbol{\beta}}\big) 
\coloneqq \lim_{t \downarrow t_{0}} \frac{H\big(\mathrm{P}_{\kern-0.1em t}^{\boldsymbol{\beta}} \, \big\vert \, \mathrm{Q}_{t}^{\boldsymbol{\beta}}\big) - H\big(\mathrm{P}_{\kern-0.1em t_{0}}^{\boldsymbol{\beta}} \, \big\vert \, \mathrm{Q}_{t_{0}}^{\boldsymbol{\beta}}\big)}{W_{2}\big(\mathrm{P}_{\kern-0.1em t}^{\boldsymbol{\beta}},\mathrm{P}_{\kern-0.1em t_{0}}^{\boldsymbol{\beta}}\big)} 
\end{equation}
along the perturbed curve $(\mathrm{P}_{\kern-0.1em t}^{\boldsymbol{\beta}})_{t_{0} \leqslant t \leqslant T}$.

\begin{theorem} \label{prop: gradient flow} Suppose \textnormal{\myassumption{assumptions.c}} hold. The following assertions hold for Lebesgue-a.e.\ $t_0 \in [0,T]$\textnormal{:} The random variables
\begin{equation} \label{eq: r.v.s.def.}
\mathcal{L}_{t_{0}}^{\uparrow} \coloneqq \nabla \log \ell_{t_{0}}^{\uparrow}(X_{t_{0}},\mathrm{P}_{\kern-0.1em t_{0}}) 
\qquad \textnormal{ and } \qquad 
\mathcal{B}_{t_{0}} \coloneqq \nabla \boldsymbol{\beta}(X_{t_0}) 
\end{equation}
are elements of $L^{2}(\mathds{P})$, and the Wasserstein metric slope of the free energy functional $\mathscr{F}$ along the McKean--Vlasov curve $(\mathrm{P}_{\kern-0.1em t})_{t_{0} \leqslant t \leqslant T}$ is given by
\begin{equation} \label{eq: w.m.s.u} 
\big\vert \partial \mathscr{F} \big\vert_{W_{2}}(\mathrm{P}_{\kern-0.1em t_{0}}) 
= - \big\Vert \mathcal{L}_{t_{0}}^{\uparrow} \big\Vert_{L^{2}(\mathds{P})}.
\end{equation}
If $\mathcal{L}_{t_{0}}^{\uparrow} + \mathcal{B}_{t_{0}} \neq 0$, the metric slope along the perturbed curve $(\mathrm{P}_{\kern-0.1em t}^{\boldsymbol{\beta}})_{t_{0} \leqslant t \leqslant T}$ is equal to
\begin{equation} \label{eq: w.m.s.p} 
\big\vert \partial \mathscr{F}^{\boldsymbol{\beta}} \big\vert_{W_{2}}\big(\mathrm{P}_{\kern-0.1em t_{0}}^{\boldsymbol{\beta}}\big) 
= - \Bigg\langle \mathcal{L}_{t_{0}}^{\uparrow} \, , \, \frac{\mathcal{L}_{t_{0}}^{\uparrow} + \mathcal{B}_{t_{0}}}{\big\Vert \mathcal{L}_{t_{0}}^{\uparrow} + \mathcal{B}_{t_{0}}\big\Vert_{L^{2}(\mathds{P})}} \Bigg\rangle_{L^{2}(\mathds{P})}.
\end{equation}
In particular, 
\begin{equation} \label{eq: inequality}
\big\vert \partial \mathscr{F} \big\vert_{W_{2}}(\mathrm{P}_{\kern-0.1em t_{0}}) 
\leqslant \big\vert \partial \mathscr{F}^{\boldsymbol{\beta}} \big\vert_{W_{2}}\big(\mathrm{P}_{\kern-0.1em t_{0}}^{\boldsymbol{\beta}}\big)
\end{equation}
with equality if and only if $\mathcal{L}_{t_{0}}^{\uparrow} + \mathcal{B}_{t_{0}}$ is a positive multiple of $\mathcal{L}_{t_{0}}^{\uparrow}$.
\begin{proof} The equality \hyperref[eq: w.m.s.u]{(\ref*{eq: w.m.s.u})} follows from \hyperref[eq: unper.rel.ent.dis.]{(\ref*{eq: unper.rel.ent.dis.})} and by taking $\boldsymbol{\beta} \equiv 0$ in \hyperref[eq: local perturbed wasserstein]{(\ref*{eq: local perturbed wasserstein})}. For the proof of \hyperref[eq: w.m.s.p]{(\ref*{eq: w.m.s.p})}, we first observe that from \hyperref[eq: per.rel.ent.dis.]{(\ref*{eq: per.rel.ent.dis.})} and \hyperref[eq: local perturbed wasserstein]{(\ref*{eq: local perturbed wasserstein})} we obtain the equality
\begin{equation} \label{eq: quo}
\big\vert \partial \mathscr{F}^{\boldsymbol{\beta}} \big\vert_{W_{2}}\big(\mathrm{P}_{\kern-0.1em t_{0}}^{\boldsymbol{\beta}}\big) 
= - \frac{   \big\Vert \mathcal{L}_{t_{0}}^{\uparrow} \big\Vert_{L^{2}(\mathds{P})}^2 +  \mathds{E}_{\mathds{P}} \Big[ \Big(  \Big\langle \nabla V + \nabla( W\ast\mathrm{P}_{\kern-0.1em t_{0}}) \, , \nabla \boldsymbol{\beta} \Big\rangle -
\Delta \boldsymbol{\beta} \Big)(X_{t_{0}})\Big] }{\big\Vert \mathcal{L}_{t_{0}}^{\uparrow} + \mathcal{B}_{t_{0}}\big\Vert_{L^{2}(\mathds{P})}}
\end{equation}
for Lebesgue-a.e.\ $t_{0} \in [0,T]$. Integrating by parts and recalling the notations in \hyperref[eq: pot]{(\ref*{eq: pot})} -- \hyperref[eq: intro.1.7]{(\ref*{eq: intro.1.7})}, we find that the expectation in the numerator of \hyperref[eq: quo]{(\ref*{eq: quo})} is equal to
\begin{equation}
\int_{\mathds{R}^{n}} \Big\langle  \log \nabla \ell_{t_{0}}^{\uparrow}(x) \, , \nabla \boldsymbol{\beta}(x) \Big\rangle  \, p_{t_0}(x) \, \mathrm{d}x 
= \Big\langle \mathcal{L}_{t_{0}}^{\uparrow} \, , \, \mathcal{B}_{t_{0}} \Big\rangle_{L^{2}(\mathds{P})}.
\end{equation}
Now \eqref{eq: inequality} follows by the Cauchy--Schwarz inequality.
\end{proof}
\end{theorem}


\subsection{A trajectorial proof of the HWBI inequality} \label{subs: hwbi}


In this subsection, we show how our trajectorial approach can be adapted to give a simple proof of the HWBI inequality. While the techniques that will be used are similar, the setting of this section is independent from the rest of the paper. 

We fix two probability measures $\nu_{0}$ and $\nu_{1}$ in $\mathscr{P}_{\mathrm{ac},2}(\mathds{R}^{n})$. By Brenier's theorem \cite{Bre91}, there exists a convex function $\varphi \colon \mathds{R}^{n} \rightarrow \mathds{R}$ such that 
\begin{equation} \label{eq: minimizer}
W^2_{2}(\nu_{0},\nu_{1}) = \int_{\mathds{R}^{n}} \vert x - \nabla \varphi(x) \vert^{2} \, \mathrm{d}\nu_{0}(x).
\end{equation}
The displacement interpolation of McCann \cite{McC97} between $\nu_{0}$ and $\nu_{1}$ is given by
\begin{equation} \label{eq: displacement interpolation}
\nu_{t}  \coloneqq  (T_t)_{\#} \nu_{0} \, , \qquad T_{t}(x) \coloneqq  (1 - t) x + t \nabla \varphi (x) \, , \qquad 0 \leqslant t \leqslant 1.
\end{equation}
In particular, since the endpoints $\nu_{0}$ and $\nu_{1}$ belong to $\mathscr{P}_{\mathrm{ac},2}(\mathds{R}^{n})$, each $\nu_{t}$ has a probability density function $\rho_{t}$; see, e.g., \cite[Remarks 5.13 (i)]{Vil03}.

As before, we consider a confinement potential $V$ and an interaction potential $W$. For each $t \in [0,1]$, we then define by analogy with \hyperref[eq: intro.1.6]{(\ref*{eq: intro.1.6})}, the $\sigma$-finite measures
\begin{equation}
\mu_{t}(A) \coloneqq \int_{A}q(x,\nu_{t})\,\mathrm{d}x 
\, , \qquad 
\mu_{t}^{\uparrow}(A) \coloneqq \int_{A}q^{\uparrow}(x,\nu_{t})\,\mathrm{d}x 
\, , 
\qquad A \in \mathscr{B}(\mathds{R}^{n}),
\end{equation}
where we recall the definitions of the density functions $q$ and $q^{\uparrow}$ in \hyperref[eq: d.f.s.]{(\ref*{eq: d.f.s.})}. In parallel to the likelihood ratio functions in \hyperref[eq: intro.1.7]{(\ref*{eq: intro.1.7})}, we define 
\begin{equation} \label{eq: r.def.}
r_{t}(x,\nu) \coloneqq \frac{\rho_{t}(x)}{q(x,\nu)} 
\, , \qquad 
r_{t}^{\uparrow}(x,\nu) \coloneqq \frac{\rho_{t}(x)}{q^{\uparrow}(x,\nu)} 
\, , \qquad
(t,x,\nu) \in [0, 1] \times \mathds{R}^{n} \times \mathscr{P}_{2}(\mathds{R}^{n}).
\end{equation}
Then the relative entropy of $\nu_{t}$ with respect to $\mu_{t}$ is given by 
\begin{equation}
H(\nu_{t} \, \vert \, \mu_{t}) 
= \int_{\mathds{R}^{n}} \rho_{t}(x) \log r_{t}(x, \nu_{t}) \,  \mathrm{d}x
\end{equation}
and the relative Fisher information of $\nu_{t}$ with respect to $\mu_{t}^{\uparrow}$ is equal to
\begin{equation}
I\big(\nu_{t} \, \vert \, \mu_{t}^{\uparrow}\big) 
= \int_{\mathds{R}^{n}} \vert \nabla \log r_{t}^{\uparrow}(x,\nu_{t})\vert^{2} \, \rho_{t}(x) \, \mathrm{d}x.
\end{equation}

We impose the following regularity conditions for \hyperref[prop: hwbi.s.]{Proposition \ref*{prop: hwbi.s.}}, noting that the strong assumptions placed on $\rho_0$ and $\rho_1$ are only temporary and will be removed in \myassumption{HWBI.assum.2} of \mythmref{thm: HWBI}.

\begin{assumptions} \label{HWBI.assum.1} The functions $V,W \colon \mathds{R}^{n} \rightarrow [0,\infty)$ are smooth and $W$ is symmetric. The probability density functions $\rho_{0}$ and $\rho_{1}$ are smooth, compactly supported and strictly positive in the interior of their respective supports. 
\end{assumptions}

\begin{proposition} \label{prop: hwbi.s.} Suppose \textnormal{\myassumption{HWBI.assum.1}} hold. Along the displacement interpolation $(\nu_{t})_{0 \leqslant t \leqslant 1}$, the rate of relative entropy dissipation at time $t=0$, with respect to the ``reference curve of probability measures'' $(\mu_{t})_{0 \leqslant t \leqslant 1}$, is given by
\begin{equation} \label{eq: hwbi.s.}
\frac{\mathrm{d}}{\mathrm{d}t} \Big\vert_{t=0}^{+} \, H(\nu_{t} \, \vert \, \mu_{t})  
= \int_{\mathds{R}^{n}} \Big\langle \nabla \log r_{0}^{\uparrow}(x,\nu_{0}) \, , \nabla \varphi(x) - x \Big\rangle \, \rho_{0}(x) \, \mathrm{d}x.
\end{equation}
\end{proposition}

Combining \hyperref[prop: hwbi.s.]{Proposition \ref*{prop: hwbi.s.}} with the displacement convexity results of McCann \cite{McC97}, we obtain the following generalization of the HWBI inequality. Equivalent versions of this inequality can be found in \cite[Theorem 4.1]{CEGH04} and \cite[Theorem D.50]{FK06}.
 
\begin{assumptions} \label{HWBI.assum.2} The functions $V,W \colon \mathds{R}^{n} \rightarrow [0,\infty)$ are smooth and $W$ is symmetric. Furthermore, $V$ and $W$ are uniformly convex, i.e., there exist real constants $\kappa_{V}$ and $\kappa_{W}$ such that
\begin{equation} \label{eq. unif.conv.}
\mathrm{Hess}(V) \geqslant \kappa_{V} I_{n} 
\, , \qquad 
\mathrm{Hess}(W) \geqslant \kappa_{W} I_{n}.
\end{equation}
\end{assumptions}

\begin{theorem} \label{thm: HWBI} Suppose \textnormal{\myassumption{HWBI.assum.2}} hold and the relative entropy $H(\nu_{1} \, \vert \, \mu_{1})$ is finite. Then
\begingroup
\addtolength{\jot}{0.7em}
\begin{align}
H(\nu_{0} \, \vert \, \mu_{0}) - H(\nu_{1} \, \vert \, \mu_{1}) \leqslant
&- \int_{\mathds{R}^{n}} \Big\langle \nabla \log r_{0}^{\uparrow}(x,\nu_{0}) \, , \nabla \varphi(x) - x \Big\rangle \, \rho_{0}(x) \, \mathrm{d}x \label{eq: HWBI} \\
&- \tfrac{\kappa_{V} + \kappa_{W}}{2} \, W_{2}^{2}(\nu_{0},\nu_{1}) + \tfrac{\kappa_{W}}{2} \, \vert b(\nu_{0}) - b(\nu_{1})\vert^{2}. \label{eq: HWBI.a}
\end{align}
\endgroup
\end{theorem} 

\begin{remark} By the Cauchy--Schwarz inequality, the right-hand side of \hyperref[eq: HWBI]{(\ref*{eq: HWBI})} can be bounded from above by
\begin{equation}
\sqrt{\int_{\mathds{R}^{n}} \big\vert \nabla \log r_{0}^{\uparrow}(x,\nu_{0})\big\vert^{2} \, \rho_{0}(x) \, \mathrm{d}x} \
\sqrt{\int_{\mathds{R}^{n}} \vert \nabla \varphi(x) - x \vert^{2} \, \rho_{0}(x) \, \mathrm{d}x} = \sqrt{I(\nu_{0} \, \vert \, \mu_{0}^{\uparrow})} \ W_{2}(\nu_{0},\nu_{1}),
\end{equation}
and we obtain the usual form of the HWBI inequality \myeqref{eq: HWBI.o}; see also \cite[Theorem 4.2]{AGK04}.
\end{remark}


\section{Proofs of the main results} \label{sec: Proofs}


This section is devoted to the proofs of the results stated in \hyperref[sec: main results]{Section \ref*{sec: main results}}. We shall first prove the main trajectorial results: \mythmref{thm: trac} and its ``perturbed'' counterpart, \mythmref{thm: p.trac}. 


\subsection{The proofs of \texorpdfstring{\hyperref[thm: trac]{Theorem \ref*{thm: trac}}}{Theorem 3.1} and \texorpdfstring{\hyperref[thm: p.trac]{Theorem \ref*{thm: p.trac}}}{Theorem 3.8}} \label{subsec: proof.traj.results}


Since \hyperref[thm: trac]{Theorem \ref*{thm: trac}} follows immediately from \hyperref[thm: p.trac]{Theorem \ref*{thm: p.trac}} by setting the perturbation $\boldsymbol{\beta} \colon \mathds{R}^{n} \rightarrow \mathds{R}$ to be the zero function, we start with the general setting of \hyperref[thm: p.trac]{Theorem \ref*{thm: p.trac}}. We first recall a classical result concerning the time reversal of diffusions.

\begin{lemma}[\textnormal{\cite[Theorem 2.1]{HP86}, \cite[Theorems G.2, G.5]{KST20b}}] \label{lem: time-reversal} Suppose \textnormal{\myassumption{assumptions.b}} hold. On $(\Omega,\mathds{G},\mathds{P}^{\boldsymbol{\beta}})$, the process 
\begin{equation} \label{eq: b.w.b.m.}
\mybar{B}_{s}^{\boldsymbol{\beta}} \coloneqq B_{T-s}^{\boldsymbol{\beta}} - B_{T}^{\boldsymbol{\beta}} - \sqrt{2} \int_{0}^{s} \nabla \log \mybar{p}_{u}^{\boldsymbol{\beta}} (\mybar{X}^{\bbeta}_{u}) \, \mathrm{d}u 
\, , \qquad 0 \leqslant s \leqslant T-t_{0}
\end{equation} 
is a Brownian motion. Moreover, the time-reversed canonical process $(\mybar{X}_{s})_{0 \leqslant s \leqslant T-t_{0}}$ satisfies 
\begin{equation} \label{eq: time-reversed diffusion}
\mathrm{d} \mybar{X}_{s} =  \Big(2 \nabla \log\mybar{\ell}_{s}^{\boldsymbol{\beta}\downarrow} - \nabla V + \nabla  (W \ast \mybar{\mathrm{P}}_{\kern-0.1em s}^{\boldsymbol{\beta}}) + \nabla \boldsymbol{\beta} \Big)(\mybar{X}_{s}) \, \mathrm{d}s 
+ \sqrt{2} \, \mathrm{d} \mybar{B}_{s}^{\boldsymbol{\beta}}.
\end{equation} 
\end{lemma}

By means of \hyperref[lem: time-reversal]{Lemma \ref*{lem: time-reversal}}, the first step in the proof of \hyperref[thm: p.trac]{Theorem \ref*{thm: p.trac}} is to compute the dynamics of the time-reversed perturbed relative entropy process \hyperref[eq: trprep]{(\ref*{eq: trprep})}. For the reader's convenience, we recall the following characterization of the L-derivative in \cite[pp.\ 383]{CD18a}.

\begin{definition} Let $f \colon \mathscr{P}_2(\R^n) \rightarrow \R$ and $\mu_0 \in \mathscr{P}_2(\R^n)$. On a probability space $(\Omega, \mathds{F}, \mathds{P})$, let $X_0$ be a random variable with distribution $\mu_0$. We define $\partial_\mu f (\mu_0) \colon \R^n \rightarrow \R^n$ as the L\emph{-derivative} of $f$ at $\mu_0$, if for any $\mu \in \mathscr{P}_2(\R^n)$ and any random variable $X$ with distribution $\mu$, 
\begin{equation*}
f(\mu) = f(\mu_0) + \mathds{E}_{\mathds{P}} \Big[\big\langle\partial_\mu f (\mu_0)(X_0) \, , \, X-X_0 \big\rangle\Big] + o\big(\Vert X- X_0 \Vert_{L^2(\mathds{P})}\big).
\end{equation*}
\end{definition}

\begin{remark} The above characterization of the L-derivative depends neither on the choice of the probability space $(\Omega, \mathds{F}, \mathds{P})$, nor of the random variables $X$ and $X_0$ used to represent $\mu$ and $\mu_0$, respectively. Moreover, if the L-derivative exists, it is uniquely defined up to $\mu_0$-equivalence. We refer to Proposition 5.25 and Remark 5.26 in \cite{CD18a} for the details.
\end{remark}

\begin{proposition} Suppose \textnormal{\myassumption{assumptions.b}} hold. On $(\Omega,\G,\P^{\bbeta})$, the time-reversed perturbed relative entropy process \textnormal{\hyperref[eq: trprep]{(\ref*{eq: trprep})}} satisfies 
\begingroup
\addtolength{\jot}{0.7em}
\begin{align}
&\mathrm{d}\log\mybar{\ell}_{s}^{\boldsymbol{\beta}}(\mybar{X}_{s}, \mybar{\mathrm{P}}_{\kern-0.1em s}^{\boldsymbol{\beta}}) 
= \Big\langle \nabla\log\mybar{\ell}_{s}^{\boldsymbol{\beta}}(\mybar{X}_{s},\mybar{\mathrm{P}}_{\kern-0.1em s}^{\boldsymbol{\beta}}) \, , \, \sqrt{2} \, \mathrm{d}\mybar{B}_{s}^{\boldsymbol{\beta}}\Big\rangle   
+ \Big( \big\vert \nabla\log\mybar{\ell}_{s}^{\boldsymbol{\beta}\downarrow} \big\vert^{2} + \tfrac{1}{2}  \big\vert \nabla( W\ast\mybar{\mathrm{P}}_{\kern-0.1em s}^{\boldsymbol{\beta}}) \big\vert^{2} \Big)(\mybar{X}_{s}) \, \mathrm{d}s 
\label{eq: ito3.p1} \\
& \qquad \qquad + \Big( \Big\langle  \tfrac{1}{2} \nabla(W\ast\mybar{\mathrm{P}}_{\kern-0.1em s}^{\boldsymbol{\beta}}) \,  , \, 2 \nabla\log\mybar{\ell}_{s}^{\boldsymbol{\beta}\downarrow} + \nabla V^{\boldsymbol{\beta}}  \Big\rangle + \langle \nabla V \, , \nabla \boldsymbol{\beta} \rangle - \Delta \boldsymbol{\beta} \Big)(\mybar{X}_{s})   \, \mathrm{d}s 
\label{eq: ito3.p3} \\
& \qquad \qquad - \mathds{E}_{\tilde{\mathds{P}}^{\boldsymbol{\beta}}}\Big[\Big\langle  \tfrac{1}{2} \nabla W(\mybar{X}_{s}-\mybar{Y}_{s})  \, , \Big(2\nabla\log\mybar{\ell}_{s}^{\boldsymbol{\beta}\downarrow} - \nabla V + \nabla ( W\ast\mybar{\mathrm{P}}_{\kern-0.1em s}^{\boldsymbol{\beta}}) + \nabla \boldsymbol{\beta} \Big)(\mybar{Y}_{s})\Big\rangle \Big] \, \mathrm{d}s, 
\label{eq: ito3.p4} 
\end{align}
\endgroup 
where $(\mybar{Y}_{s})_{0 \leqslant s \leqslant T-t_{0}}$ is a copy of the process $(\mybar{X}_{s})_{0 \leqslant s \leqslant T-t_{0}}$ on a copy $(\tilde{\Omega},\tilde{\mathds{G}},\tilde{\mathds{P}}^{\boldsymbol{\beta}})$ of the original probability space $(\Omega,\mathds{G},\mathds{P}^{\boldsymbol{\beta}})$.
\begin{proof} Applying a generalized version of It\^{o}'s formula for McKean--Vlasov diffusions \cite[Proposition 5.102]{CD18a} and using the backward dynamics in \hyperref[eq: time-reversed diffusion]{(\ref*{eq: time-reversed diffusion})}, we obtain
\begingroup
\addtolength{\jot}{0.7em}
\begin{align}
&\mathrm{d}\log\mybar{\ell}_{s}^{\boldsymbol{\beta}}(\mybar{X}_{s}, \mybar{\mathrm{P}}_{\kern-0.1em s}^{\boldsymbol{\beta}}) 
= \Big\langle \nabla\log\mybar{\ell}_{s}^{\boldsymbol{\beta}}(\mybar{X}_{s},\mybar{\mathrm{P}}_{\kern-0.1em s}^{\boldsymbol{\beta}}) \, , \, \sqrt{2} \, \mathrm{d}\mybar{B}_{s}^{\bbeta}\Big\rangle 
+ \Big(\partial_{s}\log\mybar{\ell}_{s}^{\boldsymbol{\beta}}+\Delta\log\mybar{\ell}_{s}^{\boldsymbol{\beta}}\Big)(\mybar{X}_{s},\mybar{\mathrm{P}}_{\kern-0.1em s}^{\boldsymbol{\beta}}) \, \mathrm{d}s  
\label{eq: ito1.p1} \\
& \qquad + \Big\langle \nabla\log\mybar{\ell}_{s}^{\boldsymbol{\beta}}(\mybar{X}_{s},\mybar{\mathrm{P}}_{\kern-0.1em s}^{\boldsymbol{\beta}}) \, , \Big(2\nabla\log\mybar{\ell}_{s}^{\boldsymbol{\beta}\downarrow} - \nabla V +\nabla( W\ast\mybar{\mathrm{P}}_{\kern-0.1em s}^{\boldsymbol{\beta}})+ \nabla \boldsymbol{\beta} \Big)(\mybar{X}_{s})\Big\rangle \, \mathrm{d}s
\label{eq: ito1.p2} \\
& \qquad + \mathds{E}_{\tilde{\mathds{P}}^{\boldsymbol{\beta}}}\Big[\Big\langle \Big(\partial_{\mu}\log\mybar{\ell}_{s}^{\boldsymbol{\beta}}(\mybar{X}_{s}, \mybar{\mathrm{P}}_{\kern-0.1em s}^{\boldsymbol{\beta}})\Big) \, , \Big(2\nabla\log\mybar{\ell}_{s}^{\boldsymbol{\beta}\downarrow} - \nabla V +\nabla ( W\ast\mybar{\mathrm{P}}_{\kern-0.1em s}^{\boldsymbol{\beta}}) + \nabla \boldsymbol{\beta} \Big)\Big\rangle(\mybar{Y}_{s}) \Big] \, \mathrm{d}s
\label{eq: ito1.p3} \\
& \qquad + \mathds{E}_{\tilde{\mathds{P}}^{\boldsymbol{\beta}}}\Big[\mathrm{trace}\Big(\partial_{y}\partial_{\mu}\log\mybar{\ell}_{s}^{\boldsymbol{\beta}}(\mybar{X}_{s},\mybar{\mathrm{P}}_{\kern-0.1em s}^{\boldsymbol{\beta}})(\mybar{Y}_{s})\Big)\Big] \, \mathrm{d}s, 
\label{eq: ito1.p4}
\end{align}
\endgroup
where $(\mybar{Y}_{s})_{0 \leqslant s \leqslant T-t_{0}}$ is a copy of the process $(\mybar{X}_{s})_{0 \leqslant s \leqslant T-t_{0}}$ on a copy $(\tilde{\Omega},\tilde{\mathds{G}},\tilde{\mathds{P}}^{\boldsymbol{\beta}})$ of the original probability space $(\Omega,\mathds{G},\mathds{P}^{\boldsymbol{\beta}})$. The L-derivative appearing in \hyperref[eq: ito1.p3]{(\ref*{eq: ito1.p3})} and \hyperref[eq: ito1.p4]{(\ref*{eq: ito1.p4})} is calculated to be
\begin{equation} \label{eq: d.p.1}
\Big(\partial_{\mu}\log\mybar{\ell}_{s}^{\boldsymbol{\beta}}(x,\mu)\Big)(y)  = \tfrac{1}{2} \big( \partial_{\mu}\left(W\ast\mu\right)(x) \big) (y) = - \tfrac{1}{2} \nabla W(x-y)
\end{equation}
for $(x,\mu,y) \in \mathds{R}^{n} \times \mathscr{P}_{2}(\mathds{R}^{n}) \times \mathds{R}^{n}$, see \cite[Section 5.2.2, Example 1]{CD18a} for the computation of the L-derivative of a function which is linear in the distribution variable. Consequently, we have
\begin{equation} \label{eq: d.p.2}
\mathrm{trace}\Big( \partial_{y}\partial_{\mu}\log\mybar{\ell}_{s}^{\boldsymbol{\beta}}(x,\mu)(y) \Big)  = - \tfrac{1}{2} \mathrm{trace}\big( \partial_{y}\nabla W(x-y)\big) = \tfrac{1}{2} \Delta W(x-y).
\end{equation}
Putting \hyperref[eq: d.p.1]{(\ref*{eq: d.p.1})} and \hyperref[eq: d.p.2]{(\ref*{eq: d.p.2})} into \hyperref[eq: ito1.p3]{(\ref*{eq: ito1.p3})} and \hyperref[eq: ito1.p4]{(\ref*{eq: ito1.p4})}, respectively, as well as using the identities
\begingroup
\addtolength{\jot}{0.7em}
\begin{align}
\partial_{s}\log\mybar{\ell}_{s}^{\boldsymbol{\beta}}(x,\mu) &= \partial_{s}\log\mybar{\ell}_{s}^{\boldsymbol{\beta}\downarrow}(x), \\
\nabla\log\mybar{\ell}_{s}^{\boldsymbol{\beta}}(x,\mu) &= \nabla\log\mybar{\ell}_{s}^{\boldsymbol{\beta}\downarrow}(x) + \tfrac{1}{2} \nabla(W\ast\mu)(x), \\
\Delta\log\mybar{\ell}_{s}^{\boldsymbol{\beta}}(x,\mu) &= \Delta\log\mybar{\ell}_{s}^{\boldsymbol{\beta}\downarrow}(x) + \tfrac{1}{2} \Delta(W\ast\mu)(x),
\end{align}
\endgroup
we obtain
\begingroup
\addtolength{\jot}{0.7em}
\begin{align}
&\mathrm{d}\log\mybar{\ell}_{s}^{\boldsymbol{\beta}}(\mybar{X}_{s}, \mybar{\mathrm{P}}_{\kern-0.1em s}^{\boldsymbol{\beta}}) 
= \Big\langle \nabla\log\mybar{\ell}_{s}^{\boldsymbol{\beta}}(\mybar{X}_{s},\mybar{\mathrm{P}}_{\kern-0.1em s}^{\boldsymbol{\beta}}) \, , \, \sqrt{2} \, \mathrm{d}\mybar{B}_{s}^{\boldsymbol{\beta}}\Big\rangle + \Big(\partial_{s}\log\mybar{\ell}_{s}^{\boldsymbol{\beta}\downarrow}+\Delta\log\mybar{\ell}_{s}^{\boldsymbol{\beta}\downarrow}\Big)(\mybar{X}_{s}) \, \mathrm{d}s  
\label{eq: ito2.p1} \\
& \qquad + \Big\langle \Big( \nabla\log\mybar{\ell}_{s}^{\boldsymbol{\beta}\downarrow} + \tfrac{1}{2} \nabla(W\ast\mybar{\mathrm{P}}_{\kern-0.1em s}^{\boldsymbol{\beta}}) \Big) \, , \Big(2\nabla\log\mybar{\ell}_{s}^{\boldsymbol{\beta}\downarrow} - \nabla V +\nabla( W\ast\mybar{\mathrm{P}}_{\kern-0.1em s}^{\boldsymbol{\beta}})+ \nabla \boldsymbol{\beta} \Big)\Big\rangle(\mybar{X}_{s}) \, \mathrm{d}s 
\label{eq: ito2.p2} \\
& \qquad -  \mathds{E}_{\tilde{\mathds{P}}^{\boldsymbol{\beta}}}\Big[\Big\langle  \tfrac{1}{2} \nabla W(\mybar{X}_{s}-\mybar{Y}_{s})  \, , \Big(2\nabla\log\mybar{\ell}_{s}^{\boldsymbol{\beta}\downarrow} - \nabla V +\nabla ( W\ast\mybar{\mathrm{P}}_{\kern-0.1em s}^{\boldsymbol{\beta}}) + \nabla \boldsymbol{\beta} \Big)(\mybar{Y}_{s})\Big\rangle \Big] \, \mathrm{d}s 
\label{eq: ito2.p3} \\
& \qquad + \tfrac{1}{2} \Big( \Delta(W\ast\mybar{\mathrm{P}}_{\kern-0.1em s}^{\boldsymbol{\beta}})(\mybar{X}_{s})  + \mathds{E}_{\tilde{\mathds{P}}^{\boldsymbol{\beta}}}\Big[\Delta W(\mybar{X}_{s}-\mybar{Y}_{s}) \Big] \Big) \, \mathrm{d}s. 
\label{eq: ito2.p4}
\end{align}
\endgroup
Regarding the expression of \hyperref[eq: ito2.p4]{(\ref*{eq: ito2.p4})}, we observe that $\Delta(W\ast\mybar{\mathrm{P}}_{\kern-0.1em s}^{\boldsymbol{\beta}})(\mybar{X}_{s}) = \mathds{E}_{\tilde{\mathds{P}}^{\boldsymbol{\beta}}}[\Delta W(\mybar{X}_{s}-\mybar{Y}_{s})]$. Finally, elementary computations based on \hyperref[eq: d.f.s.]{(\ref*{eq: d.f.s.})}, \hyperref[eq: perturbed kinetic FPK]{(\ref*{eq: perturbed kinetic FPK})} and \hyperref[eq: perturbed  likeilhood rato, perturbed  stationary distribution]{(\ref*{eq: perturbed  likeilhood rato, perturbed  stationary distribution})} show that the perturbed log-likelihood ratio function $(s,x) \mapsto \log\mybar{\ell}_{s}^{\boldsymbol{\beta}\downarrow}(x)$ of \textnormal{\hyperref[eq: perturbed  likeilhood rato, perturbed  stationary distribution]{(\ref*{eq: perturbed  likeilhood rato, perturbed  stationary distribution})}} satisfies
\begin{equation} \label{eq: p.b.k.e.p.}
\begin{aligned} 
\partial_{s} \log \mybar{\ell}_{s}^{\boldsymbol{\beta}\downarrow} 
&= \Big\langle \nabla\log\mybar{\ell}_{s}^{\boldsymbol{\beta}\downarrow} \, , \nabla V - \nabla( W \ast\mybar{\mathrm{P}}_{\kern-0.1em s}^{\boldsymbol{\beta}}) - \nabla \boldsymbol{\beta} \Big\rangle 
- \big\vert \nabla\log\mybar{\ell}_{s}^{\boldsymbol{\beta}\downarrow} \big\vert^{2} - \Delta\log\mybar{\ell}_{s}^{\boldsymbol{\beta}\downarrow} \\
& \qquad + \Big\langle \nabla V \, , \nabla  (W \ast\mybar{\mathrm{P}}_{\kern-0.1em s}^{\boldsymbol{\beta}}) + \nabla \boldsymbol{\beta} \Big\rangle  - \Delta( W\ast\mybar{\mathrm{P}}_{\kern-0.1em s}^{\boldsymbol{\beta}}) - \Delta \boldsymbol{\beta}
\end{aligned}
\end{equation}
on $(0,T-t_{0}) \times \mathds{R}^{n}$, with terminal condition $\log \mybar{\ell}_{T-t_{0}}^{\boldsymbol{\beta}\downarrow} = \log \mybar{\ell}_{T-t_{0}}^{\downarrow}$. Inserting \hyperref[eq: p.b.k.e.p.]{(\ref*{eq: p.b.k.e.p.})} into \hyperref[eq: ito2.p1]{(\ref*{eq: ito2.p1})}, we obtain \hyperref[eq: ito3.p1]{(\ref*{eq: ito3.p1})} -- \hyperref[eq: ito3.p4]{(\ref*{eq: ito3.p4})}. 
\end{proof}
\end{proposition}

Setting the perturbation $\boldsymbol{\beta}$ to be the zero function, we obtain the following result.

\begin{corollary} Suppose \textnormal{\myassumption{assumptions}} hold. On $(\Omega,\mathds{G},\mathds{P})$, the time-reversed relative entropy process \textnormal{\hyperref[eq: t.r.r.e.pr]{(\ref*{eq: t.r.r.e.pr})}} satisfies 
\begingroup
\addtolength{\jot}{0.7em}
\begin{align}
&\mathrm{d}\log\mybar{\ell}_{s}(\mybar{X}_{s}, \mybar{\mathrm{P}}_{\kern-0.1em s}) 
= \Big\langle \nabla\log\mybar{\ell}_{s}(\mybar{X}_{s},\mybar{\mathrm{P}}_{\kern-0.1em s}) \, , \, \sqrt{2} \, \mathrm{d}\mybar{B}_{s}\Big\rangle   
\label{eq: ito3.p1b} \\
& \qquad + \Big( \big\vert \nabla\log\mybar{\ell}_{s}^{\downarrow} \big\vert^{2} + \tfrac{1}{2}  \big\vert \nabla( W\ast\mybar{\mathrm{P}}_{\kern-0.1em s}) \big\vert^{2} 
+ \Big\langle  \tfrac{1}{2} \nabla(W\ast\mybar{\mathrm{P}}_{\kern-0.1em s}) \,  , \, 2\nabla\log\mybar{\ell}_{s}^{\downarrow}  + \nabla V  \Big\rangle \Big)(\mybar{X}_{s}) \, \mathrm{d}s 
\label{eq: ito3.p2b} \\
& \qquad -  \mathds{E}_{\tilde{\mathds{P}}}\Big[\Big\langle  \tfrac{1}{2} \nabla W(\mybar{X}_{s}-\mybar{Y}_{s})  \, , \Big(2\nabla\log\mybar{\ell}_{s}^{\downarrow}  - \nabla V +\nabla ( W\ast\mybar{\mathrm{P}}_{\kern-0.1em s})  \Big)(\mybar{Y}_{s})\Big\rangle \Big] \, \mathrm{d}s.
\label{eq: ito3.p4b} 
\end{align}
\endgroup 
Here, the process 
\begin{equation} \label{eq: b.w.b.m.sec}
\mybar{B}_{s} \coloneqq B_{T-s} - B_{T} - \sqrt{2} \int_{0}^{s} \nabla \log \mybar{p}_{u} (\mybar{X}_{u}) \, \mathrm{d}u 
\, , \qquad 0 \leqslant s \leqslant T
\end{equation} 
is a $\P$-Brownian motion with respect to the backward filtration $\G$, and $(\mybar{Y}_{s})_{0 \leqslant s \leqslant T}$ is a copy of the process $(\mybar{X}_{s})_{0 \leqslant s \leqslant T}$ on a copy $(\tilde{\Omega},\tilde{\mathds{G}},\tilde{\mathds{P}})$ of the original probability space $(\Omega,\mathds{G},\mathds{P})$.
\end{corollary}

Before turning to the final part of the proof of \mythmref{thm: trac}, we state a classical result based on the general theory of the Cameron--Martin--Maruyama--Girsanov transformation \cite{LS01}. The connection between relative entropy (the left-hand side of \hyperref[eq: ent.and.en]{(\ref*{eq: ent.and.en})} below) and energy (the right-hand side of \hyperref[eq: ent.and.en]{(\ref*{eq: ent.and.en})}) is the foundation of Föllmer’s entropy approach to the time reversal of diffusion processes on Wiener space \cite{Foe85,Foe86,Foe88}. We denote by $\mathds{W}_{\kern-0.1em x}$ the Wiener measure on $\Omega = C([0,T];\mathds{R}^{n})$ with starting point $x \in \mathds{R}^{n}$, and define by 
\begin{equation} 
\W_{\kern-0.1em x,2}(A) \coloneqq \mathds{W}_{\kern-0.1em x}\Big(\omega \in \Omega \colon  (\sqrt{2} X)(\omega) \in A\Big) \, , \qquad A \in \Bs (\Omega)
\end{equation}
the Wiener measure with starting point $x$ and variance $2$.

\begin{lemma} \label{lem: fin.ent} The relative entropy of $\mathds{P}$ with respect to $\W_{\kern-0.1em \mathrm{P}_{\kern-0.1em 0},2} \coloneqq \int_{\mathds{R}^{n}} \mathds{W}_{\kern-0.1em x,2} \, \mathrm{d}\mathrm{P}_{\kern-0.1em 0}(x)$ is given by
\begin{equation} \label{eq: ent.and.en}
H\big( \mathds{P} \, \vert \, \W_{\kern-0.1em \mathrm{P}_{\kern-0.1em 0},2}\big) 
=  \mathds{E}_{\mathds{P}}\bigg[  \int_{0}^{T} \big\vert \nabla V(X_{t}) + \nabla (W \ast \mathrm{P}_{\kern-0.1em t})(X_{t}) \big\vert^{2} \, \mathrm{d}t\bigg] < \infty.
\end{equation}
\begin{proof} Recalling \hyperref[eq: pot]{(\ref*{eq: pot})}, the drift of the McKean--Vlasov dynamics \hyperref[eq: intro.1.1]{(\ref*{eq: intro.1.1})} can be expressed as 
\begin{equation}
- \nabla \Psi^{\uparrow}(x,\mathrm{P}_{\kern-0.1em t}) = - \big( \nabla V(x) + \nabla (W \ast \mathrm{P}_{\kern-0.1em t})(x) \big) \, , \qquad  (t,x) \in [0,T] \times \mathds{R}^{n}.
\end{equation}
For any $t \in [0,T]$, using the elementary inequality $(a + b)^2 \leqslant 2 a^2 + 2b^2$, we have 
\begin{equation}  \label{eq: theta.bound}
\mathds{E}_{\mathds{P}}\Big[ \vert \nabla \Psi^{\uparrow}(X_{t},\mathrm{P}_{\kern-0.1em t}) \vert^{2}  \Big] 
\leqslant  2 \, \mathds{E}_{\mathds{P}}\Big[ \vert \nabla V(X_{t})\vert^{2} \Big] 
         + 2 \, \mathds{E}_{\mathds{P}}\Big[ \vert \nabla (W \ast \mathrm{P}_{\kern-0.1em t})(X_{t})\vert^{2} \Big].    
\end{equation}
Using the linear growth condition \hyperref[lin.gro.con.]{(\ref*{lin.gro.con.})} from \myassum{assumptions}{r.a.1}, we find 
\begin{equation} \label{fi.ent.alm.}
\mathds{E}_{\mathds{P}}\Big[ \vert \nabla V(X_{t})\vert^{2} \Big] 
\leqslant 2C^2 \Big(1 + \mathds{E}_{\mathds{P}}\big[  \vert X_{t} \vert^{2} \big] \Big) 
\leqslant 2C^2 \Big( 1 +  \mathds{E}_{\mathds{P}}\Big[ \sup_{0 \leqslant t \leqslant T} \vert X_{t} \vert^{2}\Big] \Big).
\end{equation}
Similarly, by Jensen's inequality and \hyperref[lin.gro.con.]{(\ref*{lin.gro.con.})}, we obtain
\begingroup
\addtolength{\jot}{0.7em} 
\begin{align}
\mathds{E}_{\mathds{P}} \Big[ \vert\nabla (W \ast \mathrm{P}_{\kern-0.1em t})(X_t) \vert^{2} \Big] \label{fi.ent.alm.a}
& \leqslant \int_{\mathds{R}^n \times \mathds{R}^{n}} \vert \nabla W(x -y) \vert^2 \, p_{t}(y) \, p_{t}(x) \, \mathrm{d}y \, \mathrm{d}x \\ 
& \leqslant 2 C^2 \bigg( 1 +   \int_{\mathds{R}^n \times \mathds{R}^{n}} \vert x -y \vert^2 \, p_{t}(y) \, p_{t}(x) \, \mathrm{d}y \, \mathrm{d}x \bigg)  \\
& \leqslant 2 C^2 \bigg( 1 + 2 \int_{\mathds{R}^n \times \mathds{R}^{n}} (\vert x \vert^2 + \vert y \vert^2 ) \, p_{t}(y) \, p_{t}(x) \, \mathrm{d}y \, \mathrm{d}x \bigg)  \\
& = 2 C^2 \Big( 1 + 4 \,  \mathds{E}_{\mathds{P}} \big[ \vert X_t \vert^2\big] \Big) 
\leqslant 8C^2 \Big( 1 +  \mathds{E}_{\mathds{P}}\Big[ \sup_{0 \leqslant t \leqslant T} \vert X_{t} \vert^{2}\Big]  \Big). \label{fi.ent.alm.b}
\end{align}
\endgroup
Altogether, we get
\begin{equation} \label{girs.dens.sec.zer}
\mathds{E}_{\mathds{P}}\bigg[ \int_{0}^{T} \vert  \nabla \Psi^{\uparrow}(X_{t},\mathrm{P}_{\kern-0.1em t}) \vert^{2} \, \mathrm{d}t \bigg] 
\leqslant 20 T \Big( 1 +  \mathds{E}_{\mathds{P}}\Big[ \sup_{0 \leqslant t \leqslant T} \vert X_{t} \vert^{2}\Big] \Big) < \infty,
\end{equation}
where the finiteness of this expression follows from the uniform second moment property \hyperref[fsmc]{(\ref*{fsmc})} of \hyperref[lem: well-posedness]{Lemma \ref*{lem: well-posedness}}. From \cite[Section 7.6.4]{LS01} we now conclude that $\mathds{P}$ is absolutely continuous with respect to $\W_{\kern-0.1em \mathrm{P}_{\kern-0.1em 0},2}$, and the Radon--Nikodym derivatives are given by
\begin{equation} \label{girs.dens.sec}
\frac{\mathrm{d}\mathds{P}}{\mathrm{d}\W_{\kern-0.1em \mathrm{P}_{\kern-0.1em 0},2}} \bigg\vert_{\mathcal{F}_{t}} 
= \exp \bigg( -2 \int_{0}^{t} \big\langle \nabla \Psi^{\uparrow}(X_{u},\mathrm{P}_{\kern-0.1em u}) \, , \, \sqrt{2} \, \mathrm{d}B_{u} \big\rangle + \int_{0}^{t}  \vert \nabla \Psi^{\uparrow}(X_{u},\mathrm{P}_{\kern-0.1em u}) \vert^{2} \, \mathrm{d}u \bigg) 
\, , \qquad 0 \leqslant t \leqslant T.
\end{equation}
The integrability property \hyperref[girs.dens.sec.zer]{(\ref*{girs.dens.sec.zer})} implies that the $\mathds{P}$-expectation of the stochastic integral in \hyperref[girs.dens.sec]{(\ref*{girs.dens.sec})} vanishes, and we obtain
\begin{equation} \label{foellmer.a}
H\big( \mathds{P} \, \vert \, \W_{\kern-0.1em \mathrm{P}_{\kern-0.1em 0},2}\big)  
= \mathds{E}_{\mathds{P}}\bigg[ \log \bigg( \frac{\mathrm{d}\mathds{P}}{\mathrm{d}\W_{\kern-0.1em \mathrm{P}_{\kern-0.1em 0},2}} \bigg)\bigg] 
= \mathds{E}_{\mathds{P}}\bigg[ \int_{0}^{T} \vert \nabla \Psi^{\uparrow}(X_{t},\mathrm{P}_{\kern-0.1em t}) \vert^{2} \, \mathrm{d}t \bigg] < \infty,
\end{equation}
which shows \hyperref[eq: ent.and.en]{(\ref*{eq: ent.and.en})}.
\end{proof}
\end{lemma}

Denoting $n$-dimensional Lebesgue measure by $\lambda$, we consider on $\Omega = C([0,T];\mathds{R}^{n})$ the $\sigma$-finite measure $\W_{\kern-0.1em \lambda,2} \coloneqq \int_{\mathds{R}^{n}} \mathds{W}_{\kern-0.1em x,2} \, \mathrm{d}\lambda(x)$, which is known as the law of the \emph{reversible Brownian motion} on $\mathds{R}^{n}$ with variance $2$; see \cite{Leo14a,Leo14b}. The fundamental property of reversible Brownian motion is that it is invariant under time reversal. This property can be formalized as follows. Let $R \colon \Omega \rightarrow \Omega$ be the pathwise time reversal operator on $\Omega$, given by $X_{s} \circ R = X_{T-s}$ for $s \in [0,T]$. For any measure $\mu$ on $\Omega$, we denote its time reversal by $\mybar{\mu} \coloneqq R_{\#}\mu$. Then we have the invariance property $\mybar{\W}_{\kern-0.1em \lambda,2} = \W_{\kern-0.1em \lambda,2}$. Let us also consider the probability measure $\W_{\kern-0.1em \mathrm{P}_{\kern-0.1em T},2} \coloneqq \int_{\mathds{R}^{n}} \mathds{W}_{\kern-0.1em x,2} \, \mathrm{d}\mathrm{P}_{\kern-0.1em T}(x)$ and its time reversal given by $\mybar{\W}_{\kern-0.1em \mathrm{P}_{\kern-0.1em T},2} = \int_{\mathds{R}^{n}} \mybar{\mathds{W}}_{\kern-0.1em x,2} \, \mathrm{d}\mathrm{P}_{\kern-0.1em T}(x)$. Then, as already noted in \cite[Remarks 3.7]{Foe85}, we have the following result.

\begin{lemma} \label{lemma.foe.rem} We have the relative entropy relations 
\begin{equation} \label{ea.ent.two} 
H\big( \mathds{P} \, \vert \, \W_{\kern-0.1em \lambda,2} \big)  
= H\big( \mathrm{P}_{\kern-0.1em 0} \, \vert \, \lambda \big) 
+ H\big( \mathds{P} \, \vert \, \W_{\kern-0.1em \mathrm{P}_{\kern-0.1em 0},2}\big) 
\end{equation}
and
\begin{equation} \label{ea.ent.two.b}
H\big( \mathds{P} \, \vert \, \mybar{\W}_{\kern-0.1em \lambda,2}\big) 
= H\big( \mathrm{P}_{\kern-0.1em T} \, \vert \, \lambda \big) 
+ H\big( \mathds{P} \, \vert \, \mybar{\W}_{\kern-0.1em \mathrm{P}_{\kern-0.1em T},2}\big).
\end{equation}
Furthermore, all these relative entropies are finite.
\begin{proof} For any $x \in \R^n$, we let $\P_{\kern-0.1em x} ( \, \cdot \, ) \coloneqq \P( \, \cdot \, \, \vert \, X_0 = x)$ denote (a version of) the conditional probability measure $\P$ given $X_0 = x$. By the chain rule for relative entropy \cite[Theorem 2]{Leo14b} we have
\begin{equation} \label{ea.ent.two.c} 
H\big( \P \, \vert \, \W_{\kern-0.1em \lambda,2} \big)  
= H\big( \mathrm{P}_{\kern-0.1em 0} \, \vert \, \lambda \big) 
+ \int_{\R^n} H\big( \P_{\kern-0.1em x} \, \vert \, \mathds{W}_{\kern-0.1em x,2} \big) \, \mathrm{d}\mathrm{P}_{\kern-0.1em 0}(x)
\end{equation}
and at the same time
\begin{equation} \label{ea.ent.two.d} 
H\big( \P \, \vert \, \W_{\kern-0.1em \mathrm{P}_{\kern-0.1em 0},2} \big)  
= H\big( \Pr_{\kern-0.1em 0} \, \vert \, \Pr_{\kern-0.1em 0} \big) 
+ \int_{\R^n} H\big( \P_{\kern-0.1em x} \, \vert \, \mathds{W}_{\kern-0.1em x,2} \big) \, \mathrm{d}\mathrm{P}_{\kern-0.1em 0}(x) 
= \int_{\R^n} H\big( \P_{\kern-0.1em x} \, \vert \, \mathds{W}_{\kern-0.1em x,2} \big) \, \mathrm{d}\mathrm{P}_{\kern-0.1em 0}(x),
\end{equation}
implying the first identity \myeqref{ea.ent.two}. Regarding the finite entropy assertions, we recall \myeqref{eq: free energy functional}, \myeqref{eq: intro.1.4} and observe that 
\begin{equation} \label{ea.ent} 
H\big( \mathrm{P}_{\kern-0.1em 0} \, \vert \, \lambda \big) = \int_{\mathds{R}^{n}} p_{0}(x) \log p_{0}(x) \, \mathrm{d}x = \mathcal{U}(\mathrm{P}_{\kern-0.1em 0}) \leqslant \Fs(\Pr_{\kern-0.1em 0}) < \infty,
\end{equation}
where the finiteness follows from \myassum{assumptions}{r.a.4}. Furthermore, from \mylemref{lem: fin.ent} we know that $H( \mathds{P} \, \vert \, \W_{\kern-0.1em \mathrm{P}_{\kern-0.1em 0},2}) < \infty$.

By the same arguments as above, \myeqref{ea.ent.two.b} follows again by the chain rule for relative entropy. Using the invariance property $\mybar{\W}_{\kern-0.1em \lambda,2} = \W_{\kern-0.1em \lambda,2}$, and as we already know that $H( \mathds{P} \, \vert \, \W_{\kern-0.1em \lambda,2}) < \infty$, it follows that $H( \mathds{P} \, \vert \, \mybar{\W}_{\kern-0.1em \lambda,2}) < \infty$. Let us recall now that $\Pr_{\kern-0.1em T} \in \mathscr{P}_{\kern-0.1em \mathrm{ac},2}(\mathds{R}^{n})$ by \mylemref{lem: well-posedness}. On the one hand, since $\Pr_{\kern-0.1em T}$ has finite second moment, $H( \Pr_{\kern-0.1em T} \, \vert \, \lambda)$ cannot take the value $-\infty$ as noted in \myremref{rm: rel ent}. On the other hand, the absolute continuity of $\Pr_{\kern-0.1em T}$ implies that $H( \Pr_{\kern-0.1em T} \, \vert \, \lambda)$ cannot take the value $+\infty$. Therefore, we conclude that $H( \mathds{P} \, \vert \, \mybar{\W}_{\kern-0.1em \mathrm{P}_{\kern-0.1em T},2}) < \infty$.
\end{proof}
\end{lemma}

We have assembled now all the ingredients needed for the proof of \hyperref[thm: trac]{Theorem \ref*{thm: trac}}. 

\begin{proof}[\bfseries \upshape Proof of \texorpdfstring{\hyperref[thm: trac]{Theorem \ref*{thm: trac}}}] Recalling the definition of the stochastic integral process $(\mybar{M}_{s})_{0 \leqslant s \leqslant T}$ in \hyperref[eq: martingale]{(\ref*{eq: martingale})} and of the cumulative Fisher information process $(\mybar{F}_{s})_{0 \leqslant s \leqslant T}$ in \hyperref[eq: c.f.i.]{(\ref*{eq: c.f.i.})}, we see that the stochastic differential of \hyperref[eq: ito3.p1b]{(\ref*{eq: ito3.p1b})} -- \hyperref[eq: ito3.p4b]{(\ref*{eq: ito3.p4b})} can be expressed as claimed in \hyperref[eq: d.m.d]{(\ref*{eq: d.m.d})}.

Since $(\mybar{M}_{s})_{0 \leqslant s \leqslant T}$ is a stochastic integral process, it is a continuous local martingale. In order to show that it is an $L^{2}(\mathds{P})$-bounded martingale, it suffices to show the integrability condition
\begin{equation} \label{wff.fec}
\mathds{E}_{\mathds{P}} \Big[ \big\langle \mybar{M}, \mybar{M} \big\rangle_{T} \Big] = \mathds{E}_{\mathds{P}} \bigg[  2 \int_{0}^{T} \big\vert \nabla \log \mybar{\ell}_{u}(\mybar{X}_{u},\mybar{\mathrm{P}}_{\kern-0.1em u}) \big\vert^{2} \, \mathrm{d}u \bigg] < \infty;
\end{equation}
see, e.g.\ \cite[Corollary IV.1.25]{RY99}. On $(\Omega,\mathds{G},\mathds{P})$, the time-reversed canonical process $(\mybar{X}_{s})_{0 \leqslant s \leqslant T}$ has backward dynamics
\begin{equation} \label{eq: time-reversed diffusion.unp}
\mathrm{d} \mybar{X}_{s} =  \mybar{\vartheta}_{s}(\mybar{X}_{s}) \, \mathrm{d}s + \sqrt{2} \, \mathrm{d} \mybar{B}_{s} \, , \qquad 0 \leqslant s \leqslant T
\end{equation} 
with initial distribution $\mathrm{P}_{\kern-0.1em T}$, where the drift term is given by 
\begin{equation} \label{drift.term}
\mybar{\vartheta}_{s}(x) \coloneqq 
\Big(2 \nabla \log\mybar{\ell}_{s}^{\downarrow} - \nabla V + \nabla  (W \ast \mybar{\mathrm{P}}_{\kern-0.1em s}) \Big)(x) 
= 2 \nabla \log\mybar{\ell}_{s}(x,\mybar{\mathrm{P}}_{\kern-0.1em s}) - \nabla V(x)  \in \mathds{R}^{n}
\end{equation}
for $(s,x) \in [0,T] \times \mathds{R}^{n}$. Therefore, in order to prove \hyperref[wff.fec]{(\ref*{wff.fec})}, it suffices to show the two integrability conditions 
\begin{equation} \label{fi.ent.alm.f.e.a}
\mathds{E}_{\mathds{P}}\bigg[ \int_{0}^{T} \vert \nabla V(\mybar{X}_{u})\vert^{2} \, \mathrm{d}u \bigg] < \infty
\qquad \textnormal{ and } \qquad 
\mathds{E}_{\mathds{P}}\bigg[ \int_{0}^{T} \vert \mybar{\vartheta}_{u}(\mybar{X}_{u})\vert^{2} \, \mathrm{d}u \bigg] < \infty.
\end{equation}
The first condition is a direct consequence of \hyperref[fi.ent.alm.]{(\ref*{fi.ent.alm.})}. From \cite[Lemma 2.6]{Foe85} we conclude that the expectation of the second condition is bounded by the relative entropy $H( \mathds{P} \, \vert \, \mybar{\W}_{\kern-0.1em \mathrm{P}_{\kern-0.1em T},2})$, which is finite on account of \mylemref{lemma.foe.rem}.

In order to complete the proof of \hyperref[thm: trac]{Theorem \ref*{thm: trac}}, it remains to show \hyperref[eq: fish.inf.exp]{(\ref*{eq: fish.inf.exp})}. To begin with, we take expectation with respect to $\mathds{P}$ in \hyperref[eq: c.f.i.]{(\ref*{eq: c.f.i.})} and invoke Fubini's theorem to interchange the $\mathds{P}$-expectation and the time integral. Applying once more Fubini's theorem, we swap the $\mathds{P}$-expectation with the $\tilde{\mathds{P}}$-expectation appearing in \hyperref[eq: f.i.p.]{(\ref*{eq: f.i.p.})}. Next, we recall \myassum{assumptions}{r.a.1} and use the symmetry of the interaction potential, which implies that $\nabla W(-x) = - \nabla W(x)$ for all $x \in \mathds{R}^{n}$. Furthermore, as the distribution of $\mybar{Y}_{u}$ under $\tilde{\mathds{P}}$ is the same as the distribution of $\mybar{X}_{u}$ under $\mathds{P}$, we deduce that
\begingroup
\addtolength{\jot}{0.7em} 
\begin{align}
\mathds{E}_{\mathds{P}}\big[\mybar{F}_{s}\big] 
&= \int_{0}^{s} \mathds{E}_{\mathds{P}} \Big[ \Big( \big\vert \nabla\log\mybar{\ell}_{u}^{\downarrow} \big\vert^{2} + \tfrac{1}{2}  \big\vert \nabla( W\ast\mybar{\mathrm{P}}_{\kern-0.1em u}) \big\vert^{2} + \Big\langle  \tfrac{1}{2} \nabla(W\ast\mybar{\mathrm{P}}_{\kern-0.1em u}) \,  , \, 2\nabla\log\mybar{\ell}_{u}^{\downarrow} + \nabla V  \Big\rangle \Big)(\mybar{X}_{u})\Big] \, \textnormal{d}u  \\
& \qquad + \int_{0}^{s} \mathds{E}_{\mathds{P}}\Big[\Big\langle  \tfrac{1}{2} \nabla( W\ast\mybar{\mathrm{P}}_{\kern-0.1em u})  \, , \Big(2\nabla\log\mybar{\ell}_{u}^{\downarrow} - \nabla V +\nabla ( W\ast\mybar{\mathrm{P}}_{\kern-0.1em u})  \Big)(\mybar{X}_{u})\Big\rangle \Big] \, \mathrm{d}u
\end{align}
\endgroup
for $0 \leqslant s \leqslant T$. Recalling the definitions in \hyperref[eq: pot]{(\ref*{eq: pot})} -- \hyperref[eq: intro.1.7]{(\ref*{eq: intro.1.7})}, we obtain
\begin{equation} \label{eq: ex.pe.fs.inf.a}
\mathds{E}_{\mathds{P}}\big[\mybar{F}_{s}\big] 
= \int_{0}^{s} \mathds{E}_{\mathds{P}}\Big[\big\vert\nabla\log\mybar{\ell}_{u}^{\uparrow}(\mybar{X}_{u},\mybar{\mathrm{P}}_{\kern-0.1em u})\big\vert^{2} \Big] \, \mathrm{d}u 
= \int_{0}^{s} I\big(\mybar{\mathrm{P}}_{\kern-0.1em u} \, \big\vert \, \mybar{\mathrm{Q}}_{u}^{\uparrow}\big) \, \mathrm{d}u < \infty \, , \qquad  0 \leqslant s \leqslant T,
\end{equation}
where the second equality is immediate from \hyperref[eq: intro.1.8]{(\ref*{eq: intro.1.8})}, and the finiteness of the expression in \hyperref[eq: ex.pe.fs.inf.a]{(\ref*{eq: ex.pe.fs.inf.a})} is justified as follows. Again, from \hyperref[eq: pot]{(\ref*{eq: pot})} -- \hyperref[eq: intro.1.7]{(\ref*{eq: intro.1.7})} we find
\begingroup
\addtolength{\jot}{0.7em} 
\begin{align}
\big\vert\nabla\log\mybar{\ell}_{s}^{\uparrow}(\mybar{X}_{s},\mybar{\mathrm{P}}_{\kern-0.1em s})\big\vert^{2}
&= \big\vert\nabla\log\mybar{\ell}_{s}(\mybar{X}_{s},\mybar{\mathrm{P}}_{\kern-0.1em s}) + \tfrac{1}{2} \nabla (W \ast \mybar{\mathrm{P}}_{\kern-0.1em s})(\mybar{X}_{s}) \big\vert^{2} \\
&\leqslant 2 \big\vert\nabla\log\mybar{\ell}_{s}(\mybar{X}_{s},\mybar{\mathrm{P}}_{\kern-0.1em s})\big\vert^{2} + \tfrac{1}{2} \big\vert\nabla (W \ast \mybar{\mathrm{P}}_{\kern-0.1em s})(\mybar{X}_{s}) \big\vert^{2}.
\end{align}
\endgroup
In light of \hyperref[wff.fec]{(\ref*{wff.fec})} and \hyperref[fi.ent.alm.a]{(\ref*{fi.ent.alm.a})} -- \hyperref[fi.ent.alm.b]{(\ref*{fi.ent.alm.b})}, we see that the expression in \hyperref[eq: ex.pe.fs.inf.a]{(\ref*{eq: ex.pe.fs.inf.a})} is finite, which in turn justifies a posteriori the former applications of Fubini's theorem. For a different proof in the 
\end{proof}

\begin{remark} In the above proofs of \hyperref[lem: fin.ent]{Lemmas \ref*{lem: fin.ent}}, \hyperref[lem: fin.ent]{\ref*{lemma.foe.rem}} and \hyperref[thm: trac]{Theorem \ref*{thm: trac}}, the linear growth condition \hyperref[lin.gro.con.]{(\ref*{lin.gro.con.})} from \myassum{assumptions}{r.a.1} was crucial. For a different proof in the (linear) setting without interaction (i.e., $W \equiv 0$), we refer to \cite{KST20a}, where the confinement potential $V$ (which is called $\Psi$ in \cite{KST20a}) is not necessarily of linear growth, but instead satisfies a weaker coercivity condition. The key in that setting is to apply Lemma 2.48 in \cite{KK21}.
\end{remark}

The proof of \hyperref[thm: p.trac]{Theorem \ref*{thm: p.trac}} is now an easy consequence.

\begin{proof}[\bfseries \upshape Proof of \texorpdfstring{\hyperref[thm: p.trac]{Theorem \ref*{thm: p.trac}}}] Recalling the definition of the process $(\mybar{M}_{s}^{\boldsymbol{\beta}})_{0 \leqslant s \leqslant T}$ in \hyperref[eq: p.mart]{(\ref*{eq: p.mart})} and of the perturbed cumulative Fisher information process $(\mybar{F}_{s}^{\boldsymbol{\beta}})_{0 \leqslant s \leqslant T}$ in \hyperref[eq: pcfip]{(\ref*{eq: pcfip})}, we see that the stochastic differential of \hyperref[eq: ito3.p1]{(\ref*{eq: ito3.p1})} -- \hyperref[eq: ito3.p4]{(\ref*{eq: ito3.p4})} can be expressed as claimed in \hyperref[eq: p.dmd]{(\ref*{eq: p.dmd})}.

As in the proof of \hyperref[thm: trac]{Theorem \ref*{thm: trac}} we will now argue that 
\begin{equation} \label{wff.fecb}
\mathds{E}_{\P^{\boldsymbol{\beta}}} \Big[ \big\langle \mybar{M}^{\boldsymbol{\beta}}, \mybar{M}^{\boldsymbol{\beta}} \big\rangle_{T-t_0} \Big] 
= \mathds{E}_{\P^{\boldsymbol{\beta}}} \bigg[ 2 \int_{0}^{T-t_0} \big\vert \nabla \log \mybar{\ell}_{u}^{\boldsymbol{\beta}}(\mybar{X}_{u},\mybar{\mathrm{P}}_{\kern-0.1em u}^{\boldsymbol{\beta}}) \big\vert^{2} \, \mathrm{d}u \bigg] < \infty,
\end{equation}
which will then imply that the stochastic integral process $(\mybar{M}_{s}^{\boldsymbol{\beta}})_{0 \leqslant s \leqslant T-t_0}$ is an $L^{2}(\P^{\bbeta})$-bounded martingale. To this end, we define the density $q^{\boldsymbol{\beta}}(x,\mu) \coloneqq \mathrm{e}^{-\Psi^{\boldsymbol{\beta}}(x,\mu)}$ for $(x,\mu) \in \mathds{R}^{n} \times \mathscr{P}_{2}(\mathds{R}^{n})$, and consider the ``doubly perturbed'' likelihood ratio function
\begin{equation} \label{eq: ellbb}
\ell_{t}^{\boldsymbol{\beta},\boldsymbol{\beta}}(x,\mu) \coloneqq \frac{p_{t}^{\boldsymbol{\beta}}(x)}{q^{\boldsymbol{\beta}}(x,\mu)} \, , \qquad (t,x) \in [t_0,T] \times \mathds{R}^{n}.
\end{equation}
As the \hyperref[assumptions]{Assumptions \ref*{assumptions}} are invariant under the passage from the potential $V$ to $V^{\boldsymbol{\beta}} = V + \boldsymbol{\beta}$, we can apply \hyperref[thm: trac]{Theorem \ref*{thm: trac}} to the potential $V^{\boldsymbol{\beta}}$ and obtain 
\begin{equation} \label{wff.fecb.t}
\mathds{E}_{\P^{\boldsymbol{\beta}}} \bigg[ 2 \int_{0}^{T-t_0} \big\vert \nabla \log \mybar{\ell}_{u}^{\boldsymbol{\beta},\boldsymbol{\beta}}(\mybar{X}_{u},\mybar{\mathrm{P}}_{\kern-0.1em u}^{\boldsymbol{\beta}}) \big\vert^{2} \, \mathrm{d}u \bigg] < \infty.
\end{equation}
Now, since $\ell_{t}^{\boldsymbol{\beta}}(x,\mu)  / \ell_{t}^{\boldsymbol{\beta},\boldsymbol{\beta}}(x,\mu)  = \mathrm{e}^{\boldsymbol{\beta}(x)}$, we observe that the difference
\begin{equation} 
\nabla \log \ell_{t}^{\boldsymbol{\beta}}(x,\mu) - \nabla \log \ell_{t}^{\boldsymbol{\beta},\boldsymbol{\beta}}(x,\mu) = \nabla \boldsymbol{\beta}(x)    
\end{equation}
is a bounded function. Together with \hyperref[wff.fecb.t]{(\ref*{wff.fecb.t})}, this implies \hyperref[wff.fecb]{(\ref*{wff.fecb})}.

It remains to check \hyperref[eq: fish.inf.exp.pert]{(\ref*{eq: fish.inf.exp.pert})}. A similar calculation as in the proof of \hyperref[thm: trac]{Theorem \ref*{thm: trac}} leads to the identity
\begin{equation} \label{id.ex.fis.fi}
\mathds{E}_{\mathds{P}^{\boldsymbol{\beta}}}\big[\mybar{F}_{s}^{\boldsymbol{\beta}}\big] 
= \int_{0}^{s} \mathds{E}_{\mathds{P}^{\boldsymbol{\beta}}} \Big[\big\vert\nabla\log\mybar{\ell}_{u}^{\boldsymbol{\beta}\uparrow}(\mybar{X}_{u},\mybar{\mathrm{P}}_{\kern-0.1em u}^{\boldsymbol{\beta}})\big\vert^{2} +   \Big(  \Big\langle \nabla V + \nabla( W\ast\mybar{\mathrm{P}}_{\kern-0.1em u}^{\boldsymbol{\beta}}) \, , \nabla \boldsymbol{\beta} \Big\rangle -\Delta \boldsymbol{\beta} \Big)(\mybar{X}_{u})\Big] \, \mathrm{d}u
\end{equation}
for $0 \leqslant s \leqslant T-t_{0}$. Repeating the reasoning of the previous paragraph for the function $\ell_{t}^{\boldsymbol{\beta}\uparrow}$ instead of $\ell_{t}^{\boldsymbol{\beta}}$, we find that
\begin{equation} 
\mathds{E}_{\mathds{P}^{\boldsymbol{\beta}}} \bigg[  \int_{0}^{T-t_{0}} \big\vert \nabla \log \mybar{\ell}_{u}^{\boldsymbol{\beta}\uparrow}(\mybar{X}^{\bbeta}_{u},\mybar{\mathrm{P}}_{\kern-0.1em u}^{\boldsymbol{\beta}}) \big\vert^{2} \, \mathrm{d}u \bigg] < \infty.
\end{equation}
Since the function
\begin{equation}
[0,T-t_{0}] \times \mathds{R}^{n} \ni (t,x) \longmapsto \Big\langle \nabla V + \nabla( W\ast\mybar{\mathrm{P}}_{\kern-0.1em t}^{\boldsymbol{\beta}}) \, , \nabla \boldsymbol{\beta} \Big\rangle(x) - \Delta \boldsymbol{\beta}(x)
\end{equation}
is bounded, we conclude that the quantity of \hyperref[id.ex.fis.fi]{(\ref*{id.ex.fis.fi})} is finite. Finally, recalling the definition \hyperref[eq: pe.fi.inf.id]{(\ref*{eq: pe.fi.inf.id})}, we arrive at \hyperref[eq: fish.inf.exp.pert]{(\ref*{eq: fish.inf.exp.pert})}. 
\end{proof}


\subsection{The proofs of \texorpdfstring{\hyperref[prop: hwbi.s.]{Proposition \ref*{prop: hwbi.s.}}}{Proposition 3.14} and \texorpdfstring{\hyperref[thm: HWBI]{Theorem \ref*{thm: HWBI}}}{Theorem 3.15}} \label{subsec: proofs.HWBI}


\begin{proof}[\bfseries \upshape Proof of \texorpdfstring{\hyperref[prop: hwbi.s.]{Proposition \ref*{prop: hwbi.s.}}}] The first step is to view the curve of probability density functions $(\rho_{t})_{0 \leqslant t \leqslant 1}$, corresponding to the displacement interpolation $(\nu_{t})_{0 \leqslant t \leqslant 1}$ of \hyperref[eq: displacement interpolation]{(\ref*{eq: displacement interpolation})}, as a solution to a continuity equation. Recalling the convex function $\varphi \colon \mathds{R}^{n} \rightarrow \mathds{R}$ of \hyperref[eq: minimizer]{(\ref*{eq: minimizer})}, we define $\mathsf{u}_{0} \colon \mathds{R}^{n} \to \mathds{R}$ by $\mathsf{u}_{0}(x) \coloneqq \varphi(x) - \vert x\vert^{2} \slash 2$; and for each $t \in (0,1]$, we let the function $\mathsf{u}_{t} \colon \mathds{R}^{n} \to  \mathds{R}$ be defined by the Hopf--Lax formula
\begin{equation}
\mathsf{u}_{t} (x) \coloneqq \inf_{y \in \mathds{R}^{n}} \Big( \mathsf{u}_{0}(y) + \frac{\vert x-y \vert^{2}}{2t} \Big).
\end{equation}
For all $t \in [0,1)$, we denote the gradient of $\mathsf{u}_{t}$ by $\mathsf{v}_{t} \coloneqq \nabla \mathsf{u}_{t}$. For $t = 0$, it is clear that $\mathsf{v}_{0} = \nabla \varphi - \mathrm{Id}$ is well-defined. For $t \in (0,1)$, the gradient $\mathsf{v}_{t}$ is defined Lebesgue-a.e.\ by \cite[Theorem 5.51 (i)]{Vil03}, and 
\begin{equation} \label{eq: velocity field}
\mathsf{v}_{t}(x) = \nabla \mathsf{u}_{0} \circ (T_{t})^{-1} (x) \, , \qquad \text{for all } x \in T_{t}(\mathds{R}^{n}),
\end{equation}
where $T_{t}$ is defined in \hyperref[eq: displacement interpolation]{(\ref*{eq: displacement interpolation})}. Note that the inverse of $T_{t}$ is well-defined because $T_{t}$ is injective; see \cite[Section 5.4.8]{Vil03}. From \hyperref[eq: velocity field]{(\ref*{eq: velocity field})} we see that $(\mathsf{v}_{t})_{0 \leqslant t  < 1}$ is the \emph{velocity field} associated with the trajectories $(T_{t})_{0 \leqslant t < 1}$, i.e.,
\begin{equation} \label{eq: i.e.}
T_{t}(x) = x + \int_{0}^{t} \mathsf{v}_{s}(T_{s}(x)) \, \mathrm{d}s \, , \qquad 0 \leqslant t < 1.
\end{equation} 
By \cite[Theorem 5.51 (ii)]{Vil03}, the curve of probability density functions $(\rho_{t})_{0 < t < 1}$ satisfies the continuity equation
\begin{equation} \label{eq: l.t.e}
\partial_{t}\rho_{t}(x)  + \mathrm{div}\big(\rho_{t}(x) \, \mathsf{v}_{t}(x)\big) = 0 \, , \qquad (t,x) \in(0,1) \times \mathds{R}^{n}.
\end{equation}

\smallskip

On a sufficiently rich probability space $(S, \mathcal{S}, \mathsf{P})$, we let $Z_{0} \colon S \rightarrow \mathds{R}^{n}$ be a random variable with probability distribution $\nu_{0}$. For each $0 < t \leqslant 1$, we let $Z_{t} \coloneqq T_{t}(Z_0)$. From \hyperref[eq: displacement interpolation]{(\ref*{eq: displacement interpolation})} we see that the random variable $Z_{t}$ has distribution $\nu_{t}$, and \hyperref[eq: i.e.]{(\ref*{eq: i.e.})} yields the representation
\begin{equation} \label{eq: z.e.}
Z_{t} = Z_{0} + \int_{0}^{t} \mathsf{v}_{s} (Z_{s}) \, \mathrm{d}s \, , \qquad 0 \leqslant t < 1.
\end{equation} 
In conjunction with \hyperref[eq: l.t.e]{(\ref*{eq: l.t.e})}, we deduce
\begin{equation}
\mathrm{d} \rho_{t}(Z_{t})
= \partial_{t} \rho_{t}(Z_{t}) + \big\langle \nabla \rho_{t}(Z_{t}) \, , \, \mathrm{d}Z_{t} \big\rangle 
= - \rho_{t}(Z_{t}) \, \mathrm{div}\big(\mathsf{v}_{t}(Z_{t})\big) \, \mathrm{d}t,
\end{equation} 
and thus 
\begin{equation} \label{eq: rho.dynamics}
\mathrm{d} \log\rho_{t}(Z_{t})
= - \mathrm{div}\big(\mathsf{v}_{t}(Z_{t})\big) \, \mathrm{d}t.
\end{equation}
Recalling the definition of the density function $q$ in \hyperref[eq: d.f.s.]{(\ref*{eq: d.f.s.})}, a similar argument as in \hyperref[eq: d.p.1]{(\ref*{eq: d.p.1})} shows that
\begin{equation} \label{eq: d.q.1}
\big(\partial_{\nu}\log q(x,\nu)\big)(y)  = \tfrac{1}{2} \big( \partial_{\nu}\left(W\ast\nu\right)(x) \big) (y) = - \tfrac{1}{2} \nabla W(x-y)
\end{equation}
for $(x,\nu,y) \in \mathds{R}^{n} \times \mathscr{P}_{2}(\mathds{R}^{n}) \times \mathds{R}^{n}$. Applying a generalized version of It\^{o}'s formula for McKean--Vlasov diffusions \cite[Proposition 5.102]{CD18a}, and using the dynamics \hyperref[eq: z.e.]{(\ref*{eq: z.e.})} as well as the L-derivative \hyperref[eq: d.q.1]{(\ref*{eq: d.q.1})}, we obtain
\begin{equation} \label{eq: aux.dynamics}
\mathrm{d}\log q(Z_{t},\nu_{t})
= - \Big\langle \nabla V + \tfrac{1}{2} \nabla (W \ast \nu_{t}) \, , \, \mathsf{v}_{t}\Big\rangle(Z_{t}) \, \mathrm{d}t 
+ \tfrac{1}{2} \mathds{E}_{\tilde{\mathsf{P}}} \Big[\Big\langle \nabla W(Z_{t} - \tilde{Z}_{t}) \, , \, \mathsf{v}_{t}(\tilde{Z}_{t})\Big\rangle \Big] \, \mathrm{d}t
\end{equation}
for $0 < t < 1$. Here, the process $(\tilde{Z}_{t})_{0 < t < 1}$ is defined on another probability space $(\tilde{S},\tilde{\mathcal{S}},\tilde{\mathsf{P}})$ such that the tuple $(S, \mathcal{S}, \mathsf{P},(Z_{t})_{0 < t < 1})$ is an exact copy of $( \tilde{S},\tilde{\mathcal{S}},\tilde{\mathsf{P}},(\tilde{Z}_{t})_{0 < t < 1})$. Now taking the difference between \hyperref[eq: rho.dynamics]{(\ref*{eq: rho.dynamics})} and \hyperref[eq: aux.dynamics]{(\ref*{eq: aux.dynamics})} gives the dynamics
\begingroup
\addtolength{\jot}{0.7em} 
\begin{align}
\log r_{t}(Z_{t}, \nu_{t}) - \log r_{0}(Z_{0}, \nu_{0}) 
&= \int_{0}^{t} \Big( \Big\langle \nabla V + \tfrac{1}{2} \nabla (W \ast \nu_{s}) \, , \, \mathsf{v}_{s}\Big\rangle(Z_{s}) - \mathrm{div}\big(\mathsf{v}_{s}(Z_{s})\big) \Big) \, \mathrm{d}s \label{eq: r.d.1} \\
&\qquad - \tfrac{1}{2} \int_{0}^{t} \mathds{E}_{\tilde{\mathsf{P}}} \Big[\Big\langle \nabla W(Z_{s} - \tilde{Z}_{s}) \, , \, \mathsf{v}_{s}(\tilde{Z}_{s})\Big\rangle \Big] \, \mathrm{d}s \label{eq: r.d.2}
\end{align}
\endgroup
of the relative entropy process $(\log r_{t}(Z_{t}, \nu_{t}))_{0 < t < 1}$. Next, let us make two observations. Firstly, integration by parts yields
\begin{equation} \label{eq: r.d.3}
\mathds{E}_{\mathsf{P}}\Big[ \mathrm{div}\big(\mathsf{v}_{t}(Z_{t})\big) \Big]
= - \mathds{E}_{\mathsf{P}}\Big[ \Big\langle  \nabla \log \rho_{t}(Z_{t}) \, , \, \mathsf{v}_{t}(Z_{t})  \Big\rangle \Big].
\end{equation}
Secondly, by applying Fubini's theorem, and using that $W$ is an even function as well as $(\tilde{Z}_{t})_{\#}\tilde{\mathsf{P}} = \nu_{t}$, we obtain the identity
\begin{equation} \label{eq: r.d.4}
\mathds{E}_{\mathsf{P}} \bigg[ \mathds{E}_{\tilde{\mathsf{P}}} \Big[\Big\langle \nabla W(Z_{t} - \tilde{Z}_{t}) \, , \, \mathsf{v}_{t}(\tilde{Z}_{t})\Big\rangle \Big] \bigg] 
= - \mathds{E}_{\mathsf{P}}\Big[  \Big\langle  \nabla (W \ast \nu_{t})(Z_{t}) \, , \, \mathsf{v}_{t}(Z_{t})  \Big\rangle \Big].
\end{equation}
Returning to \hyperref[eq: r.d.1]{(\ref*{eq: r.d.1})}, \hyperref[eq: r.d.2]{(\ref*{eq: r.d.2})}, we take $\mathsf{P}$-expectations and use \hyperref[eq: r.d.3]{(\ref*{eq: r.d.3})}, \hyperref[eq: r.d.4]{(\ref*{eq: r.d.4})} to obtain
\begingroup
\addtolength{\jot}{0.7em} 
\begin{align}
H(\nu_{t} \, \vert \, \mu_{t}) - H(\nu_{0} \, \vert \, \mu_{0})  
&= \int_{0}^{t} \mathds{E}_{\mathsf{P}}\Big[ \Big\langle \nabla \log \rho_{s} + \nabla V + \nabla (W \ast \nu_{s}) \, , \, \mathsf{v}_{s} \Big\rangle(Z_{s})  \Big] \, \mathrm{d}s \\
&= \int_{0}^{t} \mathds{E}_{\mathsf{P}}\Big[ \Big\langle \nabla \log r_{s}^{\uparrow}(Z_{s},\nu_{s}) \, , \, \mathsf{v}_{s}(Z_{s}) \Big\rangle  \Big] \, \mathrm{d}s,
\end{align}
\endgroup
where for the second equality we recall the notations in \hyperref[eq: r.def.]{(\ref*{eq: r.def.})} and \hyperref[eq: d.f.s.]{(\ref*{eq: d.f.s.})}. Finally, letting $t \downarrow 0$, we get
\begin{equation} 
\frac{\mathrm{d}}{\mathrm{d}t} \Big\vert_{t=0}^{+} \, H(\nu_{t} \, \vert \, \mu_{t})  
= \int_{\mathds{R}^{n}} \Big\langle \nabla \log r_{0}^{\uparrow}(x,\nu_{0}) \, , \, \mathsf{v}_{0}(x) \Big\rangle \, \rho_{0}(x) \, \mathrm{d}x;
\end{equation}
and since $\mathsf{v}_{0} = \nabla \varphi - \mathrm{Id}$, we arrive at \hyperref[eq: hwbi.s.]{(\ref*{eq: hwbi.s.})}.
\end{proof}


\begin{proof}[\bfseries \upshape Proof of \texorpdfstring{\hyperref[thm: HWBI]{Theorem \ref*{thm: HWBI}}}] Without loss of generality, we assume that the probability density functions $\rho_{0}$ and $\rho_{1}$ satisfy the strong regularity \hyperref[HWBI.assum.1]{Assumptions \ref*{HWBI.assum.1}}. The general case then follows by a density argument. We will not provide the details here, but refer to \cite[Chapter 9.4]{Vil03}, where this regularization is carried out in the simpler setting of the HWI inequality.

Let us recall the energy functionals $\mathcal{U}$, $\mathcal{V}$, $\mathcal{W}$ defined in \hyperref[eq: intro.1.4]{(\ref*{eq: intro.1.4})}, and introduce the functions
\begin{equation}
f(t) \coloneqq \mathcal{U}(\rho_{t}) \, , \qquad g(t) \coloneqq \mathcal{V}(\rho_{t}) \, , \qquad h(t) \coloneqq \mathcal{W}(\rho_{t}) \, , \qquad 0 \leqslant t \leqslant 1,
\end{equation}
where $(\rho_{t})_{0 \leqslant t \leqslant 1}$ is the curve of probability density functions corresponding to the displacement interpolation $(\nu_{t})_{0 \leqslant t \leqslant 1}$ of \hyperref[eq: displacement interpolation]{(\ref*{eq: displacement interpolation})}. Then the sum $F \coloneqq f + g + h$ of these functions satisfies the relation $F(t) = H( \nu_{t} \, \vert \, \mu_{t})$. In light of \cite[Theorem 5.15 (i)]{Vil03}, the internal energy functional $\mathcal{U}$ is displacement convex, i.e., 
\begin{equation} \label{eq: i.e.es}
f''(t) \geqslant 0 \, , \qquad 0 \leqslant t \leqslant 1.
\end{equation}
By \myassumption{HWBI.assum.2}, the confinement potential $V \colon \mathds{R}^{n} \rightarrow [0,\infty)$ is $\kappa_{V}$-uniformly convex. Therefore, \cite[Theorem 5.15 (ii)]{Vil03} implies that the potential energy functional $\mathcal{V}$ is $\kappa_{V}$-uniformly displacement convex. In other words,
\begin{equation} \label{eq: p.e.}  
g''(t) \geqslant \kappa_{V} \, W_{2}^{2}(\nu_{0},\nu_{1}) \, , \qquad 0 \leqslant t \leqslant 1.   
\end{equation}
Again from \myassumption{HWBI.assum.2}, the interaction potential $W \colon \mathds{R}^{n} \rightarrow [0,\infty)$ is assumed to be symmetric and $\kappa_{W}$-uniformly convex. Therefore, a similar argument as in the proof of \cite[Theorem 5.15 (iii)]{Vil03} leads to the $\kappa_{W}(W_{2}^{2}(\nu_{0},\nu_{1}) - \vert b(\nu_{0})-b(\nu_{1})\vert^{2})$-uniform convexity of $h$, so
\begin{equation} \label{eq: in.e.} 
h''(t) \geqslant \kappa_{W}\Big(W_{2}^{2}(\nu_{0},\nu_{1}) - \vert b(\nu_{0})-b(\nu_{1})\vert^{2}\Big) \, , \qquad 0 \leqslant t \leqslant 1. 
\end{equation}
The details of the proof of \myeqref{eq: in.e.} are postponed to \hyperref[sec: pf.in.e.]{Subsection \ref*{sec: pf.in.e.}}. By combining the estimates \hyperref[eq: i.e.es]{(\ref*{eq: i.e.es})} -- \hyperref[eq: in.e.]{(\ref*{eq: in.e.})}, we deduce that the relative entropy function $[0,1] \ni t \mapsto F(t) = H( \nu_{t} \, \vert \, \mu_{t})$ satisfies
\begin{equation} \label{eq: to.e.} 
F''(t) \geqslant (\kappa_{V}+\kappa_{W}) \, W_{2}^{2}(\nu_{0},\nu_{1}) - \kappa_{W} \, \vert b(\nu_{0})-b(\nu_{1}) \vert^{2}.
\end{equation}
Furthermore, from \hyperref[prop: hwbi.s.]{Proposition \ref*{prop: hwbi.s.}} we have
\begin{equation} \label{eq: to.e.a} 
F'(0^{+}) = \int_{\mathds{R}^{n}} \Big\langle \nabla \log r_{0}^{\uparrow}(x,\nu_{0}) \, , \nabla \varphi(x) - x \Big\rangle \, \rho_{0}(x) \, \mathrm{d}x.
\end{equation}
In conjunction with \hyperref[eq: to.e.]{(\ref*{eq: to.e.})} and \hyperref[eq: to.e.a]{(\ref*{eq: to.e.a})}, the Taylor formula  $F(1) = F(0) + F'(0^{+}) + \int_{0}^{1} (1-t) F''(t) \, \mathrm{d}t$ now yields the inequality \hyperref[eq: HWBI]{(\ref*{eq: HWBI})} -- \hyperref[eq: HWBI.a]{(\ref*{eq: HWBI.a})}.
\end{proof}


\section{Proofs of auxiliary results} \label{app}


\subsection{Proof of \texorpdfstring{\hyperref[lem: well-posedness]{Lemma \ref*{lem: well-posedness}}}{Lemma 2.4}} \label{sec: pf.well-pose}


The generalized potential $\Psi^{\uparrow}$ of \hyperref[eq: pot]{(\ref*{eq: pot})} allows us to cast the McKean--Vlasov dynamics of \hyperref[eq: intro.1.1]{(\ref*{eq: intro.1.1})} in the more compact form 
\begin{equation} \label{eq: mkv.c.f}
\mathrm{d}X_{t} = - \nabla \Psi^{\uparrow}(X_{t},\mathrm{P}_{\kern-0.1em t}) \, \mathrm{d}t + \sqrt{2} \, \mathrm{d}B_{t} \, , \qquad 0 \leqslant t \leqslant T.
\end{equation}
Then, for any two pairs $(x,\mu),(x^{\prime},\mu^{\prime}) \in \mathds{R}^{n} \times \mathscr{P}_{2}(\mathds{R}^{n})$, using the Lipschitz continuity of $\nabla V$ in \myassum{assumptions}{r.a.1} yields
\begin{equation} \label{eq: w.p.r.0}
\vert \nabla \Psi^{\uparrow}(x,\mu) - \nabla \Psi^{\uparrow} (x^{\prime},\mu^{\prime})\vert 
\leqslant \Vert \nabla V \Vert_{\mathrm{Lip}} \, \vert x - x^{\prime} \vert 
+ \vert \nabla (W\ast\mu)(x) -  \nabla (W\ast\mu^{\prime})(x^{\prime}) \vert.   
\end{equation}
For the convolution term, using Jensen's inequality and  the Lipschitz continuity of $\nabla W$ in \myassum{assumptions}{r.a.1} leads to
\begin{equation} \label{eq: w.p.r.1}
\vert \nabla (W\ast\mu)(x) -  \nabla (W\ast\mu^{\prime})(x^{\prime}) \vert
\leqslant \Vert \nabla W \Vert_{\mathrm{Lip}} \, \vert x - x^{\prime} \vert
+ \vert \nabla (W\ast\mu)(x^{\prime})  - \nabla (W\ast\mu^{\prime})(x^{\prime}) \vert.
\end{equation} 
For the last term above, by the Kantorovich--Rubinstein theorem \cite[Theorem 1.14]{Vil03}, we have
\begingroup
\addtolength{\jot}{0.7em} 
\begin{align} 
\bigg\vert \int_{\mathds{R}^{n}} \nabla W (x^{\prime} - \, \cdot \, ) \, \mathrm{d} (\mu - \mu^{\prime}) \bigg\vert  
& \leqslant \Vert \nabla W \Vert_{\mathrm{Lip}} \, \sup \bigg\{\int_{\mathds{R}^{n}} \varphi \, \mathrm{d} (\mu - \mu^{\prime}) \colon \varphi \in L^{1}(\vert \mu - \mu^{\prime} \vert),  \Vert \varphi \Vert_{\mathrm{Lip}} \leqslant 1 \bigg\} \nonumber \\
&  =        \Vert \nabla W \Vert_{\mathrm{Lip}} \, W_{1} (\mu,\mu^{\prime})
\leqslant   \Vert \nabla W \Vert_{\mathrm{Lip}} \, W_{2} (\mu,\mu^{\prime}), \label{eq: ineq}
\end{align} 
\endgroup
where
\begin{equation} 
W_{1}(\mu,\mu^{\prime}) =  \inf_{\scriptscriptstyle Y \sim \mu, Z \sim \mu^{\prime}} \mathds{E} \vert Y - Z \vert  \, , \qquad \mu, \mu^{\prime} \in \mathscr{P}_{1}(\mathds{R}^{n})
\end{equation}
denotes the $1$-Wasserstein-distance, and the inequality in \myeqref{eq: ineq} follows from Jensen's inequality. Altogether, we obtain
\begin{equation} \label{eq: w.p.r.3}
\vert \nabla \Psi^{\uparrow} (x, \mu) - \nabla \Psi^{\uparrow} (x^{\prime}, \mu^{\prime}) \vert 
\leqslant ( \Vert \nabla V \Vert_{\mathrm{Lip}} + \Vert \nabla W \Vert_{\mathrm{Lip}} ) \, \vert x - x^{\prime} \vert 
+ \Vert \nabla W \Vert_{\mathrm{Lip}} \, W_{2} (\mu, \mu^{\prime}).  
\end{equation}
In particular, this shows that the function $-\nabla \Psi^{\uparrow}$ is Lipschitz continuous on the product metric space $(\mathds{R}^{n}, \vert \, \cdot \, \vert ) \times (\mathscr{P}_{2}(\mathds{R}^{n}), W_2)$. In conjunction with \myassum{assumptions}{r.a.4}, \cite[Theorem 4.21]{CD18a} implies that the McKean--Vlasov SDE \hyperref[eq: mkv.c.f]{(\ref*{eq: mkv.c.f})} has a pathwise unique, strong solution satisfying the uniform second moment condition \hyperref[fsmc]{(\ref*{fsmc})}. Now we can linearize \myeqref{eq: mkv.c.f} by fixing the time-marginals $(\mathrm{P}_{\kern-0.1em t})_{0 \leqslant t \leqslant T}$, so that the drift term can be viewed as a function $(t, x) \mapsto \nabla \Psi^{\uparrow}(x, \mathrm{P}_{\kern-0.1em t})$, and \myeqref{eq: mkv.c.f} becomes an ordinary SDE with a time-inhomogeneous drift coefficient. 

The absolute continuity of the time-marginals $(\mathrm{P}_{\kern-0.1em t})_{0 \leqslant t \leqslant T}$ is immediate from \hyperref[lem: fin.ent]{Lemma \ref*{lem: fin.ent}}. A standard argument using the classical It\^{o}'s formula shows that the curve of probability density functions $(p_{t})_{0 \leqslant t \leqslant T}$ is a weak solution of the granular media equation \myeqref{eq: intro.1.2}. Finally, we turn to the regularity of this solution. From \myeqref{eq: w.p.r.3}, we see that the drift $x \mapsto \nabla \Psi^{\uparrow}(x, \mathrm{P}_{\kern-0.1em t})$ is Lipschitz continuous for every $t \in [0, T]$, and \myassum{assumptions}{r.a.1} implies that the drift is also of linear growth. The desired smoothness of $(p_{t})_{0 \leqslant t \leqslant T}$ now follows from a straightforward adaptation of the theorem in \cite{Rog85}, see also Remarks (i) -- (ii) therein.
\hfill \qed


\subsection{Proof of \texorpdfstring{\hyperref[eq: in.e.]{(\ref*{eq: in.e.})}}{(4.73)}} \label{sec: pf.in.e.}

We first rewrite the interaction energy functional $\mathcal{W}$ along the displacement interpolation $(\nu_{t})_{0 \leqslant t \leqslant 1}$. Using \myeqref{eq: displacement interpolation}, for any $t \in [0, 1]$, we have
\begingroup
\addtolength{\jot}{0.7em} 
\begin{align} 
h(t) 
&= \tfrac{1}{2} \int_{\R^n \times \R^n} W(x-y) \, \nu_t(\dd x) \, \nu_t(\dd y) \\
&= \tfrac{1}{2} \int_{\R^n \times \R^n} W\big(T_t(x)-T_t(y)\big) \, \nu_0(\dd x) \, \nu_0(\dd y) \\
&= \tfrac{1}{2} \int_{\R^n \times \R^n} W\Big(x-y - t\big( \theta(x) - \theta(y) \big)\Big) \, \nu_0(\dd x) \, \nu_0(\dd y),
\end{align}
\endgroup
where $\theta \colon \R^n \rightarrow \R^n$ is defined as $\theta(x) \coloneqq x - \nabla \varphi(x)$. Now, for any $t_1, t_2, \sigma \in [0, 1]$, by the $\kappa_W$-uniform convexity of $W$ in \myassumption{HWBI.assum.2}, we obtain
\begingroup
\addtolength{\jot}{0.7em} 
\begin{align} 
& \sigma h(t_1) + (1-\sigma)h(t_2) - h\big(\sigma t_1 + (1 - \sigma) t_2\big) \\
= \, & \tfrac{1}{2} \int_{\R^n \times \R^n} \bigg( \sigma W\Big(x-y - t_1\big( \theta(x) - \theta(y) \big)\Big) + (1 - \sigma) W\Big(x-y - t_2\big( \theta(x) - \theta(y) \big)\Big) \\
  & \hspace{2.2cm} -W\Big(x-y - \big(\sigma t_1 + (1 - \sigma) t_2\big)\big( \theta(x) - \theta(y) \big)\Big) \bigg) \, \nu_0(\dd x) \, \nu_0(\dd y)\\
\geqslant \, & \tfrac{1}{4} \kappa_W \sigma (1 - \sigma) (t_1 - t_2)^2 \int_{\R^n \times \R^n}  \vert \theta(x) - \theta(y) \vert^2  \, \nu_0(\dd x) \, \nu_0(\dd y). \label{eq: h.uniform.convex}
\end{align}
\endgroup
Next, we express the integral in \myeqref{eq: h.uniform.convex} as
\begingroup
\addtolength{\jot}{0.7em} 
\begin{align*}
\tfrac{1}{2} \int_{\R^n \times \R^n}  \vert \theta(x) - \theta(y) \vert^2  \, \nu_0(\dd x) \, \nu_0(\dd y) 
&= \int_{\R^n} \vert \theta(x) \vert^2 \, \nu_0(\dd x) -  \bigg\vert \int_{\R^n} \theta(x) \, \nu_0(\dd x) \bigg\vert^2 \\
&= \int_{\R^n} \vert x - \nabla \varphi(x) \vert^2 \, \nu_0(\dd x)- \bigg\vert\int_{\R^n} x \, \nu_0(\dd x) - \int_{\R^n}  x \, \nu_1(\dd x) \bigg\vert^2 \\
&= W_2^2(\nu_0, \nu_1) - \vert b(\nu_0) - b(\nu_1) \vert^2.
\end{align*}
\endgroup
Putting this back into \myeqref{eq: h.uniform.convex}, we deduce that $h$ is uniformly convex, with constant
\begin{equation}
\kappa_{W}\Big(W_{2}^{2}(\nu_{0},\nu_{1}) - \vert b(\nu_{0})-b(\nu_{1})\vert^{2}\Big).
\end{equation}

\hfill \qed


\section*{Acknowledgements} 


We are grateful to Ioannis Karatzas and Walter Schachermayer for suggesting this problem and giving us generous advice. We also thank Robert Fernholz, Miguel Garrido, Tomoyuki Ichiba, Donghan Kim, Kasper Larsen, and Mete Soner for helpful comments during the INTECH research meetings. 

B. Tschiderer acknowledges support by the Austrian Science Fund (FWF) under grant P28661, by the Vienna Science and Technology Fund (WWTF) through project MA16-021, and additionally appreciates travel support through the National Science Foundation (NSF) under grant NSF-DMS-14-05210. L.C. Yeung acknowledges support under grant NSF-DMS-20-04997.


\bibliographystyle{alpha}
{\footnotesize
\bibliography{references}}


\end{document}